\documentclass[12pt,draftcls,onecolumn]{IEEEtran}

\usepackage{amsfonts}
\usepackage{amssymb}
\usepackage{amsmath}
\usepackage{mathptmx}
\usepackage{cancel}
\usepackage{pstricks}
\usepackage{psfrag}
\usepackage{mathrsfs}
\usepackage{graphicx}
\usepackage{pgf,tikz}
\usepackage{algorithm}
\usepackage{algpseudocode}
\usepackage{soul}
\usetikzlibrary{arrows}
\def\be{\begin{equation}}
\def\ee{\end{equation}}
\def\bea{\begin{eqnarray}}
\def\eea{\end{eqnarray}}
\def\beann{\begin{eqnarray*}}
\def\eeann{\end{eqnarray*}}

\def\ns{\hspace{-1mm}}
\newcommand{\real}{{\mathbb{R}}}

\def\spacingset#1{\def\baselinestretch{#1}\small\normalsize}
\spacingset{1.35}

\newtheorem{lemma}{Lemma}
\newtheorem{theorem}{Theorem}
\newtheorem{remark}{Remark}
\newtheorem{corollary}{Corollary}

\newtheorem{problem}{Problem}
\newtheorem{example}{Example}[section]

\def\be{\begin{equation}}
\def\ee{\end{equation}}
\def\bea{\begin{eqnarray}}
\def\eea{\end{eqnarray}}
\def\beann{\begin{eqnarray*}}
\def\eeann{\end{eqnarray*}}

\def\ns{\hspace{-1mm}}

\def\proof{\noindent{\bf{\em Proof:}\ \ }}
\def\QED{\mbox{\rule[0pt]{1.5ex}{1.5ex}}}
\def\endproof{\hspace*{\fill}~\QED\par\endtrivlist\unskip}

\newcommand{\ima}{\operatorname{im}}

\newcommand{\defi}{\stackrel{\text{\tiny def}}{=}}

\newcommand{\complex}{{\mathbb{C}}}

\def\gD{{\cal D}}

\def\gH{{\cal H}}
\def\gI{{\cal I}}
\def\gJ{{\cal J}}

\def\gL{{\cal L}}
\def\gM{{\cal M}}
\def\gN{{\cal N}}

\def\gP{{\cal P}}
\def\gQ{{\cal Q}}
\def\gR{{\cal R}}
\def\gS{{\cal S}}

\def\gU{{\cal U}}
\def\gV{{\cal V}}
\def\gW{{\cal W}}
\def\gX{{\cal X}}
\def\gY{{\cal Y}}
\def\gZ{{\cal Z}}

\def\bmat{\left[ \begin{array}}
\def\emat{\end{array} \right]}

\def\bmat{\left[ \begin{array}}
\def\emat{\end{array} \right]}
\def\bsmat{\left[ \begin{smallmatrix}}
\def\esmat{\end{smallmatrix} \right]}

\def\l{{\lambda}}
\def\gP{{\cal P}}

\def\gU{{\cal U}}
\def\gL{{\cal L}}
\def\gR{{\cal R}}

\def\gV{{\cal V}}
\def\gH{{\cal H}}
\def\gS{{\cal S}}

\def\gX{{\cal X}}

\def\v{{v}}
\def\w{{w}}
\def\x{{{x}}}
\def\u{{u}}
\def\y{{y}}
\def\z{{z}}

\def\g{{g}}

\def\e{{e}}

\def\p{{p}}

\def\field{\mathbb{F}}

\newcommand{\spanR}{\operatorname{span}}

\begin{document}
\title{\LARGE{On the well-posedness in the solution of the disturbance decoupling by dynamic output feedback with self bounded and self hidden subspaces}}

\author{Fabrizio Padula, and Lorenzo Ntogramatzidis
 \thanks{Fabrizio Padula and L. Ntogramatzidis are with the Department of Mathematics and
Statistics, Curtin University, Perth,
Australia. E-mail: {\tt \{Fabrizio.Padula,L.Ntogramatzidis\}@curtin.edu.au. }}\\




}


\maketitle

\vspace{-1cm}

\IEEEpeerreviewmaketitle

\begin{abstract}
This paper studies the disturbance decoupling problem by dynamic output feedback with required closed-loop stability, in the general case of nonstrictly-proper systems. We will show that the extension of the geometric solution based on the ideas of self boundedness and self hiddenness, which is the one shown to maximize the number of assignable eigenvalues of the closed-loop, 
presents structural differences with respect to the strictly proper case. The most crucial aspect that emerges in the general case is the issue of the well-posedness of the feedback interconnection, which obviously has no counterpart in the strictly proper case. A fundamental property of the feedback interconnection that has so far remained unnoticed in the literature is investigated in this paper: 
the well-posedness condition is decoupled from the remaining solvability conditions. An important consequence of this fact is that 
the well-posedness condition written with respect to the supremal output nulling and infimal input containing subspaces does not need to be modified when we consider the solvability conditions of the problem with internal stability (where one would expect the well-posedness condition to be expressed in terms of supremal stabilizability and infimal detectability subspaces),
 and also when we consider the solution which uses the dual lattice structures of Basile and Marro.
\end{abstract}


\section{Introduction}
\label{secintro} 
The disturbance decoupling problem (DDP) played a central role in the development of the geometric approach in systems and control theory. Indeed, from the pioneering papers \cite{Basile-M-69a,Wonham-M-70}, it was recognized that geometry is a natural language for this type of problems; consequently, the 
 solvability conditions of the first disturbance decoupling problems
considered in the literature were expressed by means of inclusions involving certain subspaces.

The basic decoupling problem, consisting of the rejection of a disturbance from the output of a system by means of a static state-feedback, was solved in \cite{Basile-M-69a} and, independently,  in \cite{Wonham-M-70}, via the introduction of {\em controlled invariant subspaces}. These subspaces were then found to be powerful tools in the understanding of many system-theoretic properties of linear time-invariant (LTI) systems and in the solution of several control problems. 
The disturbance decoupling problem by static state feedback with the extra requirement of internal stability of the closed-loop was taken into account  in \cite{Wonham-M-70} with the introduction of stabilizability subspaces. An alternative solution to the same problem was suggested by Basile and Marro in \cite{Basile-M-82}, relying on the concept of self bounded controlled invariance, which, unlike the stabilizability subspaces of \cite{Wonham-M-70}, does not require  
eigenspace computation; in other words, the solution with self boundedness remains at the fundamental level of finite arithmetics.

A key contribution to the understanding of the advantages deriving from the adoption of self bounded controlled invariant subspaces in the solution of the disturbance decoupling problem by static state-feedback was given in \cite{Malabre-MG-DMC-97}, where it was shown that in the solution of this problem there is a number of closed-loop eigenvalues that are fixed for any feedback matrix which solves the decoupling problem; these unassignable eigenvalues are called the {\em fixed poles} of the decoupling problem.
It is shown in \cite{Malabre-MG-DMC-97} that choosing a particular self bounded subspace, denoted by $\gV_m$ in \cite{Basile-M-92},  is the best choice in terms of pole assignment, because it ensures that the maximum number of eigenvalues of the closed-loop can be freely assigned. 

For systems whose state is not accessible, a state-feedback decoupling filter cannot be implemented. This led to the formulation of the disturbance decoupling problem by dynamic output feedback. 
 The first paper which provided a solution to this problem is \cite{Schumacher-80-TAC-N6}. Around the same time, the same problem with the additional requirement of internal stability was addressed in \cite{Willems-C-81} and \cite{Imai-A-81}.
In \cite{Basile-MP-87-JOTAI}, an alternative geometric solution was proposed for this problem which uses self bounded subspaces, as well as their duals, the so-called {\em self hidden subspaces}. Again, the importance of this solution lies in the fact that it does not require eigenspace computation. 
Even more importantly, in \cite{DMCuellar-M-01} it was proved that this solution based on the idea of self boundedness and self hiddenness, is still the best in terms of assignability of the closed-loop dynamics, see also
\cite{DMCuellar-97} and \cite{DMCuellar-M-ECC}.

Most of the literature in geometric control has been developed for strictly proper systems, i.e., for those systems which have zero feedthrough between the input and the output. 
For a systematic and well-organized extension of the geometric approach for systems with a possibly non-zero direct feedthrough term we refer to the monograph \cite{Trentelman-SH-01}. The 
disturbance decoupling problem with dynamic output feedback and nonzero feedthrough has been completely solved in terms of stabilizability and detectability subspaces in  \cite{Stoorvogel-W-91}. More recently, the approach based on self boundedness and self hiddenness has been generalized in \cite{N-TAC-08} for the disturbance decoupling problem with static state-feedback. In \cite{N-TAC-08}, the result of \cite{Malabre-MG-DMC-97} on the fixed poles was also generalized to nonstrictly proper systems.

A significantly more challenging task is the solution of the
disturbance decoupling problem by dynamic output feedback for
nonstrictly proper systems using the concepts of self boundedness
and self hiddenness. An issue of well-posedness arises in
the case where the feedthrough between the control input
and the measurement output is non-zero. It was observed in \cite{Stoorvogel-W-91} that the solvability conditions, when dealing with the problem in its full generality, need to take into account the well-posedness: this results in a condition that cannot be expressed as the typical subspace inclusion of most control/estimation problems for which a geometric solution is available. In this paper, we study the role that the well-posedness condition plays in the disturbance decoupling problem by
dynamic output feedback.  We prove, in particular, that this condition is
invariant with respect to the stabilizing pair of self bounded and self hidden subspaces involved in the solution of the disturbance decoupling problem. In other words, we show that the well-posedness condition is disjoint, and therefore independent, from the remaining solvability conditions of the decoupling problem. This new property is the key to a full generalization of the solution of the disturbance decoupling problem by dynamic output feedback, as it shows that the fundamental requirement of stability does not reduce the set of well-posed feedback interconnections; therefore, choosing self bounded and self hidden subspaces does not impact on the solvability of the disturbance decoupling problem by dynamic output feedback. Furthermore, it also implies that the solution of \cite{Stoorvogel-W-91} can be conveniently re-written with a well-posedness condition for the supremal output nulling and infimal input containing subspaces instead of the corresponding stabilizability and detectability subspaces.

{\bf Notation.} Given a vector space $\gX$, we denote by $0_{\scriptscriptstyle \gX}$ the origin of $\gX$. The image and the kernel of matrix $A$ are denoted by $\ima\,A$ and $\ker\,A$, respectively. When $A$ is square, we denote by $\sigma(A)$ the spectrum of $A$. If $A: \gX \longrightarrow \gY$ is a linear map and if $\gJ \subseteq \gX$, the restriction of the map $A$ to $\gJ$ is denoted by $A\,|\gJ$. If $\gX=\gY$ and $\gJ$ is $A$-invariant, the eigenstructure of $A$ restricted to $\gJ$ is denoted by $\sigma\,(A\,| \gJ)$. If $\gJ_1$ and $\gJ_2$ are $A$-invariant subspaces and $\gJ_1\,{\subseteq}\,\gJ_2$, the mapping induced by $A$ on the quotient space $\gJ_2 / \gJ_1$ is denoted by $A\,| {\gJ_2}/{\gJ_1}$, and its spectrum is denoted by $\sigma\,(A\,| {\gJ_2}/{\gJ_1})$.
Given a map $A: \gX \longrightarrow \gX$ and a subspace $\gS$ of $\gX$, we denote by $\langle A\,|\, \gS \rangle$ the smallest $A$-invariant subspace of $\gX$ containing $\gS$ and by $\langle \gS\,|\, A \rangle$ the largest $A$-invariant subspace contained in $\gS$. 

\section{Problem Statements}
In what follows, whether the underlying system evolves in continuous or discrete time makes only minor differences and, accordingly, the time index set of any signal is denoted by $\mathbb{T}$, on the understanding that this represents either $\real^+$ in the continuous time or $\mathbb{N}$ in the discrete time. The symbol $\complex_g$ denotes either the open left-half complex plane $\complex^-$ in the continuous time or the open unit disc $\complex^\circ$ in the discrete time. A matrix $M \in \real^{n \times n}$ is said to be {\em asymptotically stable} if $\sigma(M) \subset \complex_g$. Finally, we say that $\l \in \complex$ is stable if $\l \in \complex_g$.
The operator
$\gD$ denotes either the time derivative in the continuous time,
i.e., $\gD\, x(t)=\dot{x}(t)$, or the unit time shift
 in the discrete time, i.e., $\gD\, x(t)=x(t+1)$. 

We consider the system $\Sigma$ governed by
\beann
\Sigma: \left\{\begin{array}{rcl}
\gD\,{\x}(t) \ns&\ns = \ns&\ns A\,\x(t)+B\,\u(t)+H\,\w(t) \\
\y(t) \ns&\ns = \ns&\ns C\,\x(t)+D_y\,\u(t)+G_y\,\w(t)\\
\z(t) \ns&\ns = \ns&\ns E\,\x(t)+D_z\,\u(t)+G_z\,\w(t),
\end{array}\right.
\eeann
where, for all $t\in \mathbb{T}$, the vector $\x(t)\in \gX=\real^n$ denotes the state, $\u(t)\in \gU=\real^m$ is the control input, $\w(t)\in \gW=\real^q$ is the disturbance input, $\y(t)\in \gY=\real^p$ is the measurement output
and $\z(t)\in \gZ=\real^r$ is the to-be-controlled output.
%
We consider also the regulator $\Sigma_{\scriptscriptstyle C}$ ruled by
\beann
\Sigma_{\scriptscriptstyle C}: \left\{\begin{array}{rcl}
\gD{\p}(t) \ns&\ns = \ns&\ns A_c\,\p(t)+B_c\,\y(t) \\
\u(t) \ns&\ns = \ns&\ns C_c\,\p(t)+D_c\,\y(t),
\end{array}\right.
\eeann
where, for all $t\in \mathbb{T}$, the vector $\p(t)\in \gP=\real^s$ is the state of the regulator. 
We want to control the system $\Sigma$ with the regulator $\Sigma_{\scriptscriptstyle C}$ such that in the closed-loop system the output $\z$
does not depend on the disturbance input $\w$. 

We say that the feedback interconnection of system $\Sigma$ with the regulator $\Sigma_{\scriptscriptstyle C}$ is {\em well posed} if the matrix $I-D_y\,D_c$ is non-singular, see \cite[Chpt.~3]{Trentelman-SH-01}. In such case, the closed-loop system $\Sigma_{\scriptscriptstyle CL}$ can be written in state-space form as
\bea
\label{cl}
\Sigma_{\scriptscriptstyle CL}: \left\{\begin{array}{rcl}
\gD{\hat{x}}(t) 
\ns&\ns = \ns&\ns \widehat{A}\,\hat{x}(t) +
\widehat{H}\, \w(t) \\
\z(t)
\ns&\ns = \ns&\ns \widehat{C}\,\hat{x}(t)+ \widehat{G}\,\w(t),
\end{array}\right.
\eea
where $\hat{x}(t)=\bsmat x(t) \\[1mm] p(t) \esmat$ is the extended state, and the matrices in (\ref{cl}) are defined by
{
\beann
\widehat{A} \ns&\ns \defi \ns&\ns \bmat{cc}
A+B D_c W\,C & B \,C_c+B\,D_c\,W\,D_y\,C_c \\
B_c\,W\,C & A_c+B_c\,W\,D_y\,C_c  \emat, \quad \widehat{H} \defi\bmat{c}  H+B\,D_c\,W\,G_y   \\
  B_c\,W\,G_y  \emat, \\
\widehat{C} \ns&\ns \defi \ns&\ns[\begin{array}{cc}
E+D_z\,D_c\,
W\,
C & D_z\,C_c+D_z\,D_c\,W\,D_y\,C_c \end{array}], \quad \widehat{G} \defi G_z+D_z\,D_c\,W\,G_y,
\eeann}
\ \\[-0.5cm]
where $W=(I-D_y\,D_c)^{-1}$.

\begin{figure}
\hspace{5cm}
\begin{tikzpicture}[line cap=round,line join=round,>=triangle 45,x=0.25cm,y=0.25cm]

\draw [color=black,ultra thick] (0,10) -- (6,10) -- (6,14) -- (0,14) -- (0,10);
\draw[color=black] (3,12) node {\Large $\Sigma$};
\draw [->,color=black,thick] (-2,11) -- (0,11);
\draw [-,color=black,thick] (6,11) -- (8,11);
\draw [-,color=black,thick] (8,11) -- (8,6);
\draw[color=black] (8.8,9) node {\large $\y$};
\draw [->,color=black,thick] (8,6) -- (6,6);
\draw [-,color=black,thick] (0,6) -- (-2,6);
\draw [-,color=black,thick] (-2,6) -- (-2,11);
\draw [color=black,ultra thick] (0,4) -- (6,4) -- (6,8) -- (0,8) -- (0,4);
\draw[color=black] (3,6) node {\Large $\Sigma_c$};
\draw[color=black] (-2.8,9) node {\large $\u$};
\draw [->,color=black,thick] (-6,13) -- (0,13);
\draw[color=black] (-6,13.7) node {\large $\w$};
\draw [->,color=black,thick] (6,13) -- (12,13);
\draw[color=black] (12,13.7) node {\large $\z$};
\draw [color=black,thin, dashed] (-4,3) -- (10,3) -- (10,15) -- (-4,15) -- (-4,3);
\draw[color=black] (12,3) node {\large $\Sigma_{\scriptscriptstyle CL}$};
\end{tikzpicture}
\end{figure}

The transfer function of the closed-loop system $\Sigma_{\scriptscriptstyle CL}$ is
\[
G_{\z,\w}(\l)=\widehat{C}\,(\l\,I-\widehat{A})^{-1}\widehat{H}+\widehat{G},
\]
where $\l$ represents the $s$ variable of the Laplace transform in the continuous time or the $z$ variable of the $\gZ$-transform in the discrete time.

In this paper we are concerned with two problems:

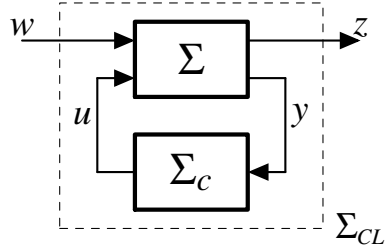
\begin{problem}{\sc [DDP by Dynamic Output  Feedback]}\\
\label{pro1}
Find a compensator $\Sigma_{\scriptscriptstyle C}$ for $\Sigma$ such that the feedback interconnection of $\Sigma$ with $\Sigma_{\scriptscriptstyle C}$ is well posed and the transfer function matrix $G_{\z,\w}(\l)$
 of the closed-loop system $\Sigma_{\scriptscriptstyle CL}$ is zero.
\end{problem}

\ \\[-.8cm]

\begin{problem}{\sc [DDP by Dynamic Output Feedback with Stability]}\\
\label{pro2}
Find a compensator $\Sigma_{\scriptscriptstyle C}$ for $\Sigma$ such that the feedback interconnection of $\Sigma$ with $\Sigma_{\scriptscriptstyle C}$ is well posed,  the transfer function matrix $G_{\z,\w}(\l)$ of the closed-loop system $\Sigma_{\scriptscriptstyle CL}$ is
 zero and all the eigenvalues of $\widehat{A}$ are in $\complex_g$.
\end{problem}

\section{Geometric background}
Consider a quadruple $(A,B,C,D)$ associated with the non-strictly proper state-space (continuous or discrete-time) system
\beann
\left\{\begin{array}{rcl}
\gD{\x}(t) \ns&\ns = \ns&\ns A\,\x(t)+B\,\u(t) \\
\y(t) \ns&\ns = \ns&\ns C\,\x(t)+D\,\u(t)\end{array}\right.
\eeann
We denote by $\gR$ the reachable subspace of the pair $(A,B)$, which is the smallest $A$-invariant subspace containing the column-space of $B$, i.e., $\gR=\langle A\,|\,\ima B\rangle$. We denote by $\gQ$ the unobservable subspace of the pair $(C,A)$, which is the largest $A$-invariant subspace contained in the null-space of $C$, i.e., $\gQ=\langle \ker C\,|\,A\rangle$.
A subspace $\gV$ is said to be an $(A,B)$-controlled invariant subspace if,
for any initial state $\x_0\in \gV$, there exists a control function $\u$ such that the state trajectory generated by the system remains identically on $\gV$; equivalently, $\gV$ is $(A,B)$-controlled invariant if the subspace inclusion  $A\,\gV\subseteq \gV+\ima B$ holds. The control function that maintains the trajectory on $\gV$ can always be expressed as a static state feedback $\u(t)=F\,\x(t)$. The condition of  $(A,B)$-controlled invariance can be equivalently expressed by saying that there exists a feedback matrix $F$ such that $(A+B\,F)\,\gV\subseteq \gV$. In this case, we say that $F$ is a {\em controlled invariant friend} of $\gV$. 
A subspace $\gV$ is said to be an $(A,B,C,D)$-output nulling subspace if,
for any initial state $\x_0\in \gV$, there exists a control function $\u$ such that the state trajectory generated by the system remains in $\gV$ and the output remains identically at zero; equivalently, $\gV$ is $(A,B,C,D)$-output nulling if the subspace inclusion  
\[
\bmat{c} A \\
C\emat \,\gV\subseteq (\gV\oplus 0_{\scriptscriptstyle \scriptscriptstyle \gY})+\ima \bmat{c}B\\ D \emat
\]
 holds. The control function that maintains the trajectory on $\gV$ can again be expressed as the static state feedback $\u(t)=F\,\x(t)$. The condition of $(A,B,C,D)$-output nullingness can be equivalently expressed by saying that there exists a feedback matrix $F$ such that 
 \[
 \bmat{c} A+B\,F \\ C+D\,F \emat \,\gV\subseteq \gV\oplus 0_{\scriptscriptstyle \scriptscriptstyle \gY}.
 \]
  In this case, we say that $F$ is an {\em output nulling friend} of $\gV$. 
  It is easy to see that if $F$ is an $(A,B,C,D)$-output nulling friend of $\gV$, we have also the inclusion (in the complexification of $\gX$)
  \bea
\label{eqV}
\ker\bigl((C+D\,F)\,(\l\,I-A-B\,F)^{-1}\bigr)\supseteq \gV
\eea
  for all $\lambda\in \complex$, see \cite{Stoorvogel-W-91}. We denote by $\mathfrak{F}_{\scriptscriptstyle (A,B,C,D)}(\gV)$ the set of $(A,B,C,D)$-output nulling friends of $\gV$. 

It is easy to see that the set of $(A,B,C,D)$-output nulling subspaces is closed under addition. Thus, we can define the largest $(A,B,C,D)$-output nulling subspace $\gV^\star_{\scriptscriptstyle (A,B,C,D)}$ (also referred to as the {\em weakly unobservable subspace}), which is the set of all initial states for which a control function exists that maintains the output identically at zero. The sequence of subspaces $(\gV_i)_{i\in \mathbb{N}}$ given by
\beann
\left\{\begin{array}{rclccccc}
\gV_0 \ns&\ns = \ns&\ns \gX \\
\gV_{i+1} \ns&\ns = \ns&\ns \bmat{c} A\\C \emat^{-1}\left((\gV_i\oplus 0_{\scriptscriptstyle \scriptscriptstyle \gY})+\ima \bmat{c} B \\ D \emat\right)
\end{array}\right.
\eeann
is monotonically non-increasing and converges to $\gV^\star_{\scriptscriptstyle (A,B,C,D)}$ in at most $n-1$ steps, i.e., 
$\gV_0\supset \gV_1\supset \ldots \supset \gV_h = \gV_{h+1}=\ldots$ implies $\gV^\star_{\scriptscriptstyle (A,B,C,D)}=\gV_h$, with $h \le n-1$.

Given an $(A,B,C,D)$-output nulling subspace $\gV$, we can define the $(A,B,C,D)$-reachability subspace $\gR_{\gV}$ on $\gV$ as the set of points that can be reached from the origin by means of control functions that maintain the trajectory on $\gV$ and the output at zero. Given an output nulling friend $F$ of $\gV$, we can determine $\gR_{\gV}$ as
\[
\gR_{\gV}=\langle A+B\,F\,|\,\gV\cap B\,\ker D\rangle.
\]
The eigenvalues of $A+B\,F$, for $F$ that varies in $\mathfrak{F}_{\scriptscriptstyle (A,B,C,D)}(\gV)$, can be divided into two multi-sets: the eigenvalues of the mapping $A+B\,F\,|\,\gV$ and the eigenvalues of $A+B\,F\,|\,\frac{\gX}{\gV}$. In turn, the eigenvalues of $A+B\,F\,|\,\gV$ can be divided into two multi-sets: the eigenvalues of $A+B\,F\,|\,\gR_{\gV}$ are all freely assignable with a suitable choice of $F\in \mathfrak{F}_{\scriptscriptstyle (A,B,C,D)}(\gV)$, whereas the eigenvalues of $A+B\,F\,|\,\frac{\gV}{\gR_{\gV}}$ are independent from $F\in \mathfrak{F}_{\scriptscriptstyle (A,B,C,D)}(\gV)$. Likewise, the eigenvalues of $A+B\,F\,|\,\frac{\gV+\gR}{\gV}$ are all freely assignable with a suitable choice of $F\in \mathfrak{F}_{\scriptscriptstyle (A,B,C,D)}(\gV)$, whereas the eigenvalues of $A+B\,F\,|\,\frac{\gX}{\gV+\gR}$ are fixed for all $F\in \mathfrak{F}_{\scriptscriptstyle (A,B,C,D)}(\gV)$. 
The {\em fixed poles} of $\gV$ can be defined as the unassignable eigenvalues of $A+B\,F$ with $F\in \mathfrak{F}_{\scriptscriptstyle (A,B,C,D)}(\gV)$, i.e.,
\beann
\label{sfixV_II}
\sigma_{\rm fixed}(\gV) \defi \sigma \left(A+B\,F\,\Big|\,\frac{\gX}{\gV+\gR}\right) \uplus 
\sigma \left(A+B\,F\,\Big|\,\frac{\gV}{\gR_{\gV}}\right). 
\eeann
It is easy to see that $\sigma_{\rm fixed}(\gV)$ can be alternatively characterized as
\bea
\label{sfixV}
\sigma_{\rm fixed}(\gV) = \sigma \left(A+B\,F\,\Big|\,\frac{\gX}{\gR}\right) \uplus 
\sigma \left(A+B\,F\,\Big|\,\frac{\gV\cap \gR}{\gR_{\gV}}\right), \qquad F\in \mathfrak{F}_{\scriptscriptstyle (A,B,C,D)}(\gV),
\eea
see \cite{DMCuellar-M-01,Malabre-MG-DMC-97}.
We say that $\gV$ is 
\begin{itemize}
\item {\em internally stabilizable} if there exists $F\in \mathfrak{F}_{\scriptscriptstyle (A,B,C,D)}(\gV)$ such that
 $\sigma(A+B\,F\,|\,\gV)\subset \complex_g$, or, equivalently,  if $\sigma(A+B\,F\,|\,\frac{\gV}{\gR_{\gV}})\subset \complex_g$; 
 \item {\em externally stabilizable} if there exists $F\in \mathfrak{F}_{\scriptscriptstyle (A,B,C,D)}(\gV)$ such that $\sigma(A+B\,F\,|\,\frac{\gX}{\gV})\subset \complex_g$, or, equivalently, if $\sigma\bigl(A+B\,F\,|\,\frac{\gX}{\gV+\gR}\bigr)\subset \complex_g$. 
\end{itemize}
An $(A,B,C,D)$-output nulling subspace that is internally stabilizable is also referred to as an $(A,B,C,D)$-{\em stabilizability output nulling subspace}: specifically, an $(A,B,C,D)$-output nulling subspace $\gV$ is an $(A,B,C,D)$-stabilizability output nulling subspace if there exists 
$F\in \mathfrak{F}_{\scriptscriptstyle (A,B,C,D)}(\gV)$ such that $\sigma(A+B\,F\,|\,\gV)\subset \complex_g$.
 The set of $(A,B,C,D)$-stabilizability output nulling subspaces is closed under addition, and thus it admits a maximum, that we denote by $\gV^\star_{\scriptscriptstyle (A,B,C,D),g}$: this subspace can be interpreted as the set of all initial states for which an input function exists that maintains the output at zero and the state trajectory converges to the origin.

An $(A,B,C,D)$-output nulling subspace $\gR$ for which an output nulling friend $F$ exists such that the spectrum of $A+B\,F\,|\,\gR$ is arbitrary is called an $(A,B,C,D)$-{\em reachability output nulling subspace}. The set of $(A,B,C,D)$-reachability output nulling subspaces is closed under addition, and thus it admits a maximum, that we denote by $\gR^\star_{\scriptscriptstyle (A,B,C,D)}$: there holds
\[
\gR^\star_{\scriptscriptstyle (A,B,C,D)}\subseteq \gV^\star_{\scriptscriptstyle (A,B,C,D),g}\subseteq \gV^\star_{\scriptscriptstyle (A,B,C,D)}.
\]
The subspace $\gR^\star_{\scriptscriptstyle (A,B,C,D)}$ is also the output nulling reachability subspace on $\gV^\star_{\scriptscriptstyle (A,B,C,D)}$, i.e., $\gR^\star_{\scriptscriptstyle (A,B,C,D)} = \gR_{\gV^\star_{\scriptscriptstyle (A,B,C,D)}}$. This subspace can be interpreted as the set of all initial states that are reachable from the origin by control inputs that maintain the output at zero. The spectrum of $A+B\,F\,|\,\frac{\gV^\star_{\scriptscriptstyle (A,B,C,D)}}{\gR^\star_{\scriptscriptstyle (A,B,C,D)}}$ is the {\em invariant zero structure} of the system, and it is denoted by $Z_{\scriptscriptstyle (A,B,C,D)}$. 

We say that an $(A,B,C,D)$-output nulling subspace $\gV$ is $(A,B,C,D)$-{\em self bounded} if, for any initial state $\x_0\in \gV$, any control that gives an identically zero output is such that the entire state trajectory is forced to evolve on $\gV$. In terms of subspace inclusions, $\gV$ is $(A,B,C,D)$-self bounded if one of the following equivalent conditions holds:
\begin{enumerate}
\item $\gV\supseteq \gV^\star_{\scriptscriptstyle (A,B,C,D)}\cap B\,\ker D$;
\item $\gV \supseteq \gR^\star_{\scriptscriptstyle (A,B,C,D)}$.
\end{enumerate}
It follows immediately that $\gR^\star_{\scriptscriptstyle (A,B,C,D)}$ and $\gV^\star_{\scriptscriptstyle (A,B,C,D)}$ are $(A,B,C,D)$-self bounded subspaces.
If $\gV_1$ and $\gV_2$ are $(A,B,C,D)$-self bounded subspaces and $\gV_1\subseteq \gV_2$, then every $(A,B,C,D)$-output nulling friend of $\gV_2$ is also an $(A,B,C,D)$-output nulling friend of $\gV_1$, i.e., $\mathfrak{F}_{\scriptscriptstyle (A,B,C,D)}(\gV_2)\subseteq 
\mathfrak{F}_{\scriptscriptstyle (A,B,C,D)}(\gV_1)$. In particular, since $\gR^\star_{\scriptscriptstyle (A,B,C,D)}\subseteq \gV^\star_{\scriptscriptstyle (A,B,C,D)}$, every $(A,B,C,D)$-output nulling friend of $\gV^\star_{\scriptscriptstyle (A,B,C,D)}$ is also an $(A,B,C,D)$-output nulling friend of $\gR^\star_{\scriptscriptstyle (A,B,C,D)}$.

Moreover, the intersection of $(A,B,C,D)$-self bounded subspaces is $(A,B,C,D)$-self bounded. Thus, if we define $\Phi_{\scriptscriptstyle (A,B,C,D)}$ to be the set of $(A,B,C,D)$-self bounded subspaces, then $\Phi_{\scriptscriptstyle (A,B,C,D)}$ admits both a maximum, which is $\gV^\star_{\scriptscriptstyle (A,B,C,D)}$, and a minimum, which is $\gR^\star_{\scriptscriptstyle (A,B,C,D)}$.\\

Most of the results on conditioned invariance are introduced by duality. We recall that the dual of a quadruple $(A,B,C,D)$ is the quadruple $(A^\top,C^\top,B^\top,D^\top)$.
A subspace $\gS$ is said to be a $(C,A)$-conditioned invariant subspace if the subspace inclusion  $A\,(\gS\cap \ker C)\subseteq \gS$ holds. 
The  $(C,A)$-conditioned invariance condition can be equivalently expressed by saying that there exists an output-injection matrix $G$ such that $(A+G\,C)\,\gS\subseteq \gS$. In this case, we say that $G$ is a {\em conditioned invariant friend} of $\gS$. 
A subspace $\gL$ is $(C,A)$-conditioned invariant subspace if and only if $\gL^\perp$ is $(A^\top,C^\top)$-controlled invariant.
A subspace $\gS$ is said to be an $(A,B,C,D)$-input containing subspace if the subspace inclusion  
\[
[\begin{array}{cc} A & B \end{array}]\left((\gS \oplus \gU) \cap \ker [\begin{array}{cc} C & D \end{array}]\right)\subseteq \gS
\]
 holds.  A subspace $\gL$ is $(A,B,C,D)$-input containing if and only if $\gL^\perp$ is $(A^\top,C^\top,B^\top,D^\top)$-output nulling.
 The condition of input containingingness can be equivalently expressed by saying that there exists an output-injection matrix $G$ such that 
 \[
 \bmat{c} A+G\,C \\ B+G\,D \emat \,(\gS\oplus \gU) \subseteq \gS.
 \]
  In this case, we say that $G$ is an $(A,B,C,D)$-{\em input containing friend} of $\gS$.
    It is easy to see that if $G$ is an $(A,B,C,D)$-input containing friend of $\gS$, we have also
\bea
\label{eqS}
\ima \bigl((\l\,I-A-G\,C)^{-1}(B+G\,D)\bigr)\subseteq \gS
\eea
  for all $\lambda\in \complex$ in the complexification of $\gX$, see \cite{Stoorvogel-W-91}. 
     We denote by $\mathfrak{G}_{\scriptscriptstyle (A,B,C,D)}(\gS)$ the set of $(A,B,C,D)$-input containing friends of $\gS$. 
The set of $(A,B,C,D)$-input containing subspaces is closed under intersection. Thus, we can define the smallest $(A,B,C,D)$-input containing subspace $\gS^\star_{\scriptscriptstyle (A,B,C,D)}$ (also referred to as the {\em strongly controllable subspace}). 
 The sequence of subspaces $(\gS_i)_{i\in \mathbb{N}}$ given by
\beann
\left\{\begin{array}{rclccccc}
\gS_0 \ns&\ns = \ns&\ns 0_{\scriptscriptstyle \scriptscriptstyle \gX} \\
\gS_{i+1} \ns&\ns = \ns&\ns  [\begin{array}{cc} A& B \end{array}]\left(({\gS}_i \oplus \gU)\cap\ker [\begin{array}{cc} C & D \end{array}]\right)
\end{array}\right.
\eeann
is monotonically non-decreasing and converges to $\gS^\star_{\scriptscriptstyle (A,B,C,D)}$ in at most $n-1$ steps, i.e., 
$\gS_0\subset \gS_1\subset \ldots \subset \gS_h = \gS_{h+1}=\ldots$ implies $\gS^\star_{\scriptscriptstyle (A,B,C,D)}=\gS_h$, with $h \le n-1$.
There holds also $\gS^\star_{\scriptscriptstyle (A,B,C,D)}=\bigl(\gV^\star_{\scriptscriptstyle (A^\top,C^\top,B^\top,D^\top)}\bigr)^\perp$.

Given an $(A,B,C,D)$-input containing subspace $\gS$ and a corresponding $(A,B,C,D)$-input containing friend $G$, we define the $(A,B,C,D)$-detectability subspace associated to it as
\[
\gQ_{\gS}=\langle \gS+C^{-1}\ima D \,|\,A+G\,C
\rangle,
\]
and is the orthogonal complement of the reachability subspace on $\gS^\perp$.
The eigenvalues of $A+G\,C$, for $G\in \mathfrak{G}_{\scriptscriptstyle (A,B,C,D)}(\gS)$, can be divided into  the eigenvalues of the mapping $A+G\,C\,|\,\gS$ and the eigenvalues of $A+G\,C\,|\,\frac{\gX}{\gS}$. In turn, the eigenvalues of $A+G\,C\,|\,\gS$ can be divided into two multi-sets: the eigenvalues of $A+G\,C\,|\,(\gS\cap \gQ)$ are fixed, whereas the eigenvalues of $A+G\,C\,|\,\frac{\gS}{\gS\cap \gQ}$
all freely assignable with a suitable choice of $G\in \mathfrak{G}_{\scriptscriptstyle (A,B,C,D)}(\gS)$. Likewise, the eigenvalues of $A+G\,C\,|\,\frac{\gQ_{\gS}}{\gS}$ are fixed, while the eigenvalues of $A+G\,C\,|\,\frac{\gX}{\gQ_{\gS}}$ are
 freely assignable with a suitable choice of $G\in \mathfrak{G}_{\scriptscriptstyle (A,B,C,D)}(\gS)$. 
 The fixed poles of $\gS$ are defined as the unassignable eigenvalues of $A+G\,C$ with $G\in \mathfrak{G}_{\scriptscriptstyle (A,B,C,D)}(\gS)$, i.e.,
\beann
\sigma_{\rm fixed}(\gS) \defi \sigma \left(A+G\,C\,\Big|\,\frac{\gQ_{\gS}}{\gS}\right) \uplus 
\sigma \left(A+G\,C\,\Big|\,\gS\cap \gQ \right),
\eeann
or, which is the same, as
\beann
\sigma_{\rm fixed}(\gS) = \sigma \left(A+G\,C\,\Big|\,\frac{\gQ_{\gS}}{\gS+\gQ}\right) \uplus 
\sigma \left(A+G\,C\,\Big|\,\gQ\right), \qquad G\in \mathfrak{G}_{\scriptscriptstyle (A,B,C,D)}(\gS).
\eeann

 We say that the $(A,B,C,D)$-input containing subspace $\gS$ is 
\begin{itemize}
\item {\em internally detectable} if there exists $G\in \mathfrak{G}_{\scriptscriptstyle (A,B,C,D)}(\gS)$ such that
 $\sigma(A+G\,C\,|\,\gS)\subset \complex_g$, or, equivalently,  if $\sigma(A+G\,C\,|\,\gS\cap \gQ)\subset \complex_g$; 
 \item {\em externally detectable} if there exists $G\in \mathfrak{G}_{\scriptscriptstyle (A,B,C,D)}(\gS)$ such that $\sigma(A+G\,C\,|\,\frac{\gX}{\gS})\subset \complex_g$, or, equivalently, if $\sigma(A+G\,C\,|\,\frac{\gQ_{\gS}}{\gS})\subset \complex_g$. 
\end{itemize}

An $(A,B,C,D)$-input containing subspace that is externally detectable is also referred to as an  $(A,B,C,D)$-{\em detectability input containing subspace}: specifically, an $(A,B,C,D)$-input containing subspace $\gS$ is an $(A,B,C,D)$-detectability input containing subspace if there exists 
$G\in \mathfrak{G}_{\scriptscriptstyle (A,B,C,D)}(\gS)$ such that $\sigma(A+G\,C\,|\,\frac{\gX}{\gS})\subset \complex_g$.

 The set of $(A,B,C,D)$-detectability input containing subspaces admits a minimum, that we denote by $\gS^\star_{\scriptscriptstyle (A,B,C,D),g}$: there holds $\gS^\star_{\scriptscriptstyle (A,B,C,D),g}=\bigl(\gV^\star_{\scriptscriptstyle (A^\top,C^\top,B^\top,D^\top),g}\bigr)^\perp$.
 
An input containing subspace $\gQ$ for which an $(A,B,C,D)$-input containing friend $G$ exists such that the spectrum of $A+G\,C\,|\,\frac{\gX}{\gQ}$ is arbitrary is called an $(A,B,C,D)$-{\em unobservability input containing subspace}. The set of $(A,B,C,D)$-unobservability input containing subspaces is closed under intersection, and thus it admits a minimum, that we denote by $\gQ^\star_{\scriptscriptstyle (A,B,C,D)}$: there holds
\[
\gS^\star_{\scriptscriptstyle (A,B,C,D)}\subseteq \gS^\star_{\scriptscriptstyle (A,B,C,D),g}\subseteq \gQ^\star_{\scriptscriptstyle (A,B,C,D)}.
\]
There holds also $\gQ^\star_{\scriptscriptstyle (A,B,C,D)} = \gQ_{\gS^\star_{\scriptscriptstyle (A,B,C,D)}}$.  The spectrum $A+G\,C\,|\,\frac{\gQ^\star_{\scriptscriptstyle (A,B,C,D)}}{\gS^\star_{\scriptscriptstyle (A,B,C,D)}}$ coincides with the {\em invariant zero structure} of the system, so that
$Z_{\scriptscriptstyle (A,B,C,D)}=\sigma(A+B\,F\,|\,\frac{\gV^\star_{\scriptscriptstyle (A,B,C,D)}}{\gR^\star_{\scriptscriptstyle (A,B,C,D)}})=\sigma(A+G\,C\,|\,\frac{\gQ^\star_{\scriptscriptstyle (A,B,C,D)}}{\gS^\star_{\scriptscriptstyle (A,B,C,D)}})$. Finally, we recall that $\gQ^\star_{\scriptscriptstyle (A,B,C,D)}$ is the dual of $\gR^\star_{\scriptscriptstyle (A,B,C,D)}$, i.e., $\gQ^\star_{\scriptscriptstyle (A,B,C,D)}=\bigl(\gR^\star_{\scriptscriptstyle (A^\top,C^\top,B^\top,D^\top)}\bigr)^\perp$. 

\begin{figure}
\begin{tikzpicture}[line cap=round,line join=round,>=triangle 45,x=0.39cm,y=0.39cm]

\draw [color=black, thin] (5,0) circle [radius=0.5];
\draw[color=black] (6.5,0) node {$0_{\scriptscriptstyle \scriptscriptstyle \gX}$};
\draw [-,color=black,thick] (5,0.5) -- (5,3.5);
\draw[color=gray, thin] (2+0.5,2) node {$\text{\sc free}  \left\{\begin{array}{c}  \phantom{\Big|}  \\  \phantom{\Big|} \end{array}\right. $};
\draw [color=black, thin] (5,4) circle [radius=0.5];
\draw[color=black] (6.5,4) node {\small $\gR_{\gV}$};
\draw [-,color=black,thick] (5,4.5) -- (5,9.5);
\draw[color=gray, thin] (1.9+0.5,7) node {$\text{\sc fixed}  \left\{\begin{array}{c} \phantom{\Big|} \\  \phantom{\Big|}  \\  \phantom{\Big|} \end{array}\right. $};
\draw [color=black, thin] (5,10) circle [radius=0.5];
\draw[color=black] (6.5,10) node {\small $\gV$};
\draw [-,color=black,thick] (5,10.5) -- (5,15.5);
\draw[color=gray, thin] (2+0.5,13) node {$\text{\sc free}  \left\{\begin{array}{c} \phantom{\Big|} \\  \phantom{\Big|}  \\  \phantom{\Big|} \end{array}\right. $};
\draw [color=black, thin] (5,16) circle [radius=0.5];
\draw[color=black] (7.5,16) node {$\gV+\gR$};
\draw [-,color=black,thick] (5,16.5) -- (5,19.5);
\draw[color=gray, thin] (1.9+0.5,18) node {$\text{\sc fixed}  \left\{\begin{array}{c}  \phantom{\Big|}  \\  \phantom{\Big|} \end{array}\right. $};
\draw [color=black, thin] (5,20) circle [radius=0.5];
\draw[color=black] (6.5,20) node {$\gX$};

\draw [color=black, thin] (5+12,0) circle [radius=0.5];
\draw[color=black] (6.5+12,0) node {$0_{\scriptscriptstyle \scriptscriptstyle \gX}$};
\draw [-,color=black,thick] (5+12,0.5) -- (5+12,3.5);
\draw[color=gray, thin] (2+0.5+12,2) node {$\text{\sc free}  \left\{\begin{array}{c}  \phantom{\Big|}  \\  \phantom{\Big|} \end{array}\right. $};
\draw [color=black, thin] (5+12,4) circle [radius=0.5];
\draw[color=black] (6.5+12+0,4) node {\small $\gR^\star$};
\draw [-,color=black,thick] (5+12,4.5) -- (5+12,9.5);
\draw [color=black, thin, fill=white] (5+12,7) circle [radius=0.5];
\draw[color=black] (6.5+12+0.2,7) node {\small $\gV^\star_{g}$};
\draw[color=gray, thin] (1.9+0.5+12,7) node {$\text{\sc fixed}  \left\{\begin{array}{c} \phantom{\Big|} \\  \phantom{\Big|}  \\  \phantom{\Big|} \end{array}\right. $};
\draw[color=gray] (0.6+12,6) node {$Z_{\scriptscriptstyle (A,B,C,D)}$};
\draw [color=black, thin, fill=white] (5+12,10) circle [radius=0.5];
\draw[color=black] (6.5+12+0,10) node {\small $\gV^\star$};
\draw [-,color=black,thick] (5+12,10.5) -- (5+12,15.5);
\draw[color=gray, thin] (2+0.5+12,13) node {$\text{\sc free}  \left\{\begin{array}{c} \phantom{\Big|} \\  \phantom{\Big|}  \\  \phantom{\Big|} \end{array}\right. $};
\draw [color=black, thin] (5+12,16) circle [radius=0.5];
\draw[color=black] (7.5+12+0.3,16) node {$\gV^\star+\gR$};
\draw [-,color=black,thick] (5+12,16.5) -- (5+12,19.5);
\draw[color=gray, thin] (1.9+0.5+12,18) node {$\text{\sc fixed}  \left\{\begin{array}{c}  \phantom{\Big|}  \\  \phantom{\Big|} \end{array}\right. $};

\draw [color=black, thin] (5+12,20) circle [radius=0.5];
\draw[color=black] (6.5+12,20) node {$\gX$};

\draw [color=black, thin] (5+24,0) circle [radius=0.5];
\draw[color=black] (6.5+24,0) node {$0_{\scriptscriptstyle \scriptscriptstyle \gX}$};
\draw [-,color=black,thick] (5+24,0.5) -- (5+24,3.5);
\draw[color=gray, thin] (1.9+0.5+24,2) node {$\text{\sc fixed}  \left\{\begin{array}{c}  \phantom{\Big|}  \\  \phantom{\Big|} \end{array}\right. $};
\draw [color=black, thin] (5+24,4) circle [radius=0.5];
\draw[color=black] (7.5+24,4) node {\small $\gS\cap \gQ$};
\draw [-,color=black,thick] (5+24,4.5) -- (5+24,9.5);
\draw[color=gray, thin] (2+0.5+24,7) node {$\text{\sc free}  \left\{\begin{array}{c} \phantom{\Big|} \\  \phantom{\Big|}  \\  \phantom{\Big|} \end{array}\right. $};
\draw [color=black, thin] (5+24,10) circle [radius=0.5];
\draw[color=black] (6.5+24,10) node {\small $\gS$};
\draw [-,color=black,thick] (5+24,10.5) -- (5+24,15.5);
\draw[color=gray, thin] (1.9+0.5+24,13) node {$\text{\sc fixed}  \left\{\begin{array}{c} \phantom{\Big|} \\  \phantom{\Big|}  \\  \phantom{\Big|} \end{array}\right. $};
\draw [color=black, thin] (5+24,16) circle [radius=0.5];
\draw[color=black] (6.5+24,16) node {$\gQ_{\gS}$};
\draw [-,color=black,thick] (5+24,16.5) -- (5+24,19.5);
\draw[color=gray, thin] (2+0.5+24,18) node {$\text{\sc free}  \left\{\begin{array}{c}  \phantom{\Big|}  \\  \phantom{\Big|} \end{array}\right. $};
\draw [color=black, thin] (5+24,20) circle [radius=0.5];
\draw[color=black] (6.5+24,20) node {$\gX$};

\draw [color=black, thin] (5+12+24,0) circle [radius=0.5];
\draw[color=black] (6.5+12+24,0) node {$0_{\scriptscriptstyle \scriptscriptstyle \gX}$};
\draw [-,color=black,thick] (5+12+24,0.5) -- (5+12+24,3.5);
\draw[color=gray, thin] (1.9+0.5+12+24,2) node {$\text{\sc fixed}  \left\{\begin{array}{c}  \phantom{\Big|}  \\  \phantom{\Big|} \end{array}\right. $};
\draw [color=black, thin] (5+12+24,4) circle [radius=0.5];
\draw[color=black] (5.5+12+2+24,4) node {\small $\gS^\star\cap \gQ$};
\draw [-,color=black,thick] (5+12+24,4.5) -- (5+12+24,9.5);
\draw [color=black, thin, fill=white] (5+12+24,13) circle [radius=0.5];
\draw[color=black] (6.5+12+0.2+24,13) node {\small $\gS^\star_{g}$};
\draw[color=gray, thin] (2+0.5+12+24,7) node {$\text{\sc free}  \left\{\begin{array}{c} \phantom{\Big|} \\  \phantom{\Big|}  \\  \phantom{\Big|} \end{array}\right. $};
\draw[color=gray] (0.6+12+24,12) node {$Z_{\scriptscriptstyle (A,B,C,D)}$};
\draw [color=black, thin, fill=white] (5+12+24,10) circle [radius=0.5];
\draw[color=black] (6.5+12+0+24,10) node {\small $\gS^\star$};
\draw [-,color=black,thick] (5+12+24,10.5) -- (5+12+24,15.5);
\draw[color=gray, thin] (1.9+0.5+12+24,13) node {$\text{\sc fixed}  \left\{\begin{array}{c} \phantom{\Big|} \\  \phantom{\Big|}  \\  \phantom{\Big|} \end{array}\right. $};
\draw [color=black, thin] (5+12+24,16) circle [radius=0.5];
\draw[color=black] (8.5+12-2+24,16) node {$\gQ^\star$};
\draw [-,color=black,thick] (5+12+24,16.5) -- (5+12+24,19.5);
\draw[color=gray, thin] (2+0.5+12+24,18) node {$\text{\sc free}  \left\{\begin{array}{c}  \phantom{\Big|}  \\  \phantom{\Big|} \end{array}\right. $};

\draw [color=black, thin] (5+12+24,20) circle [radius=0.5];
\draw[color=black] (6.5+12+24,20) node {$\gX$};
\draw [color=black, thin, fill=white] (5+12+24,13) circle [radius=0.5];

\end{tikzpicture}
\label{lattici}
\caption{Hasse diagrams of output nulling and input containing subspaces. All the subspaces refer to the quadruple $(A,B,C,D)$, but the pedix $_{\scriptscriptstyle (A,B,C,D)}$ has been dropped for the sake of readability.}
\end{figure}
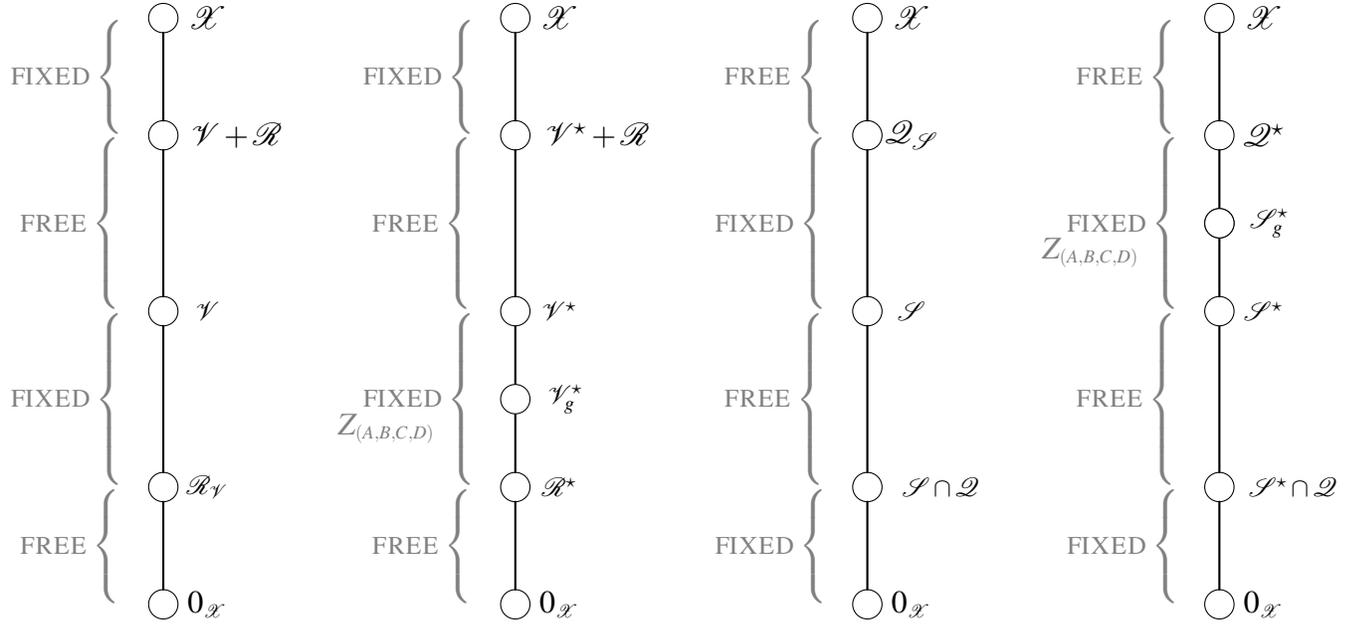

We say that an $(A,B,C,D)$-input containing subspace $\gS$ is $(A,B,C,D)$-{\em self hidden} if one of the following equivalent conditions holds:
\begin{enumerate}
\item $\gS \subseteq \gS^\star_{\scriptscriptstyle (A,B,C,D)}+C^{-1} \ima  D$;
\item $\gS \subseteq \gQ^\star_{\scriptscriptstyle (A,B,C,D)}$.
\end{enumerate}
Thus, $\gQ^\star_{\scriptscriptstyle (A,B,C,D)}$ and $\gS^\star_{\scriptscriptstyle (A,B,C,D)}$ are $(A,B,C,D)$-self hidden subspaces.
If $\gS_1$ and $\gS_2$ are $(A,B,C,D)$-self hidden subspaces and $\gS_1\subseteq \gS_2$, then every $(A,B,C,D)$-input containing friend of $\gS_1$ is also an $(A,B,C,D)$-input containing friend of $\gS_2$, i.e., $\mathfrak{G}_{\scriptscriptstyle (A,B,C,D)}(\gS_1)\subseteq 
\mathfrak{G}_{\scriptscriptstyle (A,B,C,D)}(\gS_2)$. In particular, every $(A,B,C,D)$-input containing friend of $\gS^\star_{\scriptscriptstyle (A,B,C,D)}$ is also an $(A,B,C,D)$-input containing friend of $\gQ^\star_{\scriptscriptstyle (A,B,C,D)}$.

Moreover, the sum of $(A,B,C,D)$-self hidden subspaces is $(A,B,C,D)$-self hidden. Thus, if we define $\Psi_{\scriptscriptstyle (A,B,C,D)}$ to be the set of $(A,B,C,D)$-self hidden subspaces, then $\Psi_{\scriptscriptstyle (A,B,C,D)}$ admits both a maximum, which is $\gQ^\star_{\scriptscriptstyle (A,B,C,D)}$, and a minimum, which is $\gS^\star_{\scriptscriptstyle (A,B,C,D)}$.\\

We recall also the two well-known identities
\beann
\gR^\star_{\scriptscriptstyle (A,B,C,D)} &=& \gV^\star_{\scriptscriptstyle (A,B,C,D)}\cap \gS^\star_{\scriptscriptstyle (A,B,C,D)}, \\
\gQ^\star_{\scriptscriptstyle (A,B,C,D)} &=& \gV^\star_{\scriptscriptstyle (A,B,C,D)}+ \gS^\star_{\scriptscriptstyle (A,B,C,D)}.
\eeann

\section{Dual lattice structures}
The following results extend the classic results that relate the concepts of output nullingness and input containingness, see  \cite[Chpt.~5]{Basile-M-92}.  

\begin{lemma}
\label{lem41}
Let $\gV$ be an $(A,B,C,D)$-output nulling subspace and let $\gS$ be an $(A,B,C,D)$-input containing subspace. Then, $\gS \supseteq B\,\ker D$ and
$\gV\subseteq C^{-1}\,\ima D$.
\end{lemma}
\proof We have
\[
B\,\ker D = [\begin{array}{cc} A & B \end{array}]\left((0_{\scriptscriptstyle \gX} \oplus \gU) \cap \ker [\begin{array}{cc} C & D \end{array}]\right)\subseteq [\begin{array}{cc} A & B \end{array}]\left((\gS\oplus \gU) \cap \ker [\begin{array}{cc} C & D \end{array}]\right)\subseteq \gS,
\]
 which proves the first.
The second can be proved by duality. 
\endproof

\begin{theorem}
\label{th2}
Let $\gV$ be an $(A,B,C,D)$-output nulling subspace and let $\gS$ be an $(A,B,C,D)$-input containing subspace. Then:
\begin{itemize}
\item $\gV\cap \gS$ is an $(A,B,C,D)$-output nulling subspace;
\item $\gV+ \gS$ is an $(A,B,C,D)$-input containing subspace.
\end{itemize}
\end{theorem}
\proof
We prove the first. Let us consider $\x\in \gV \cap \gS$. Since $\x\in \gV$, there exist $\x_v\in \gX$ and ${\omega}\in \gU$ such that
$\bsmat A\\[1mm] C\esmat \,\x=\bsmat \x_v \\[1mm] 0 \esmat+\bsmat B \\[1mm] D \esmat \,{\omega}$,
which can be written as the two equations
\bea
A\,\x &= &\x_v +B\,{\omega}, \label{o1} \\
C\,\x &=&\phantom{\x_v +}D\,{\omega}. \label{o2} 
\eea
 Since $\x\in \gS$, there exist $\x_s\in \gX$ and $\u\in \gU$ such that
$[\begin{array}{cc} A & B \end{array}]\,\bsmat  \x \\[1mm] \u \esmat=\x_s$
and $C\,\x+D\,\u={0}$,
which can be written as
\bea
A\,\x +B\,\u & = & \x_s, \label{o3} \\
C\,\x +D\,\u&= & 0. \label{o4} 
\eea
Subtracting (\ref{o1}) to (\ref{o3}) gives
$\x_s-\x_v=B\,({\omega}-\u)$, and subtracting (\ref{o2}) to (\ref{o4}) gives
$D\,(\omega-\u)={0}$, so that
$\x_s-\x_v\in B\,\ker D\subseteq \gS$.
 It follows that $\x_v\in \gS$. From (\ref{o1}-\ref{o2}), it follows that 
$\bsmat A\\[1mm] C\esmat \,\x \in \left((\gV\cap \gS) \oplus 0_{\scriptscriptstyle \gY}\right) +\ima \bsmat B\\[1mm] D\esmat$,
 and since $\x\in \gV\cap \gS$, the subspace $\gV\cap \gS$ is $(A,B,C,D)$-output nulling. 
 The second can be proved by duality.
\endproof

%
We now consider the two quadruples $(A,B,E,D_z)$ and $(A,[\begin{array}{cc} B & H \end{array}],E,[\begin{array}{cc} D_z & G_z \end{array}])$. We denote by $(\hat{\gV}_i)_{i \in \mathbb{N}}$ and $(\tilde{\gV}_i)_{i \in \mathbb{N}}$ the two sequences that converge in at most $n-1$ steps to $\gV^\star_{\scriptscriptstyle (A,B,E,D_z)}$ and $\gV^\star_{\scriptscriptstyle (A,[\,B \;\; H \,],E,[\, D_z \;\; G_z \,])}$, respectively. Similarly, we denote by $\gS^\star_{\scriptscriptstyle (A,B,E,D_z)}$ and $\gS^\star_{\scriptscriptstyle (A,[\,B \;\; H \,],E,[\, D_z \;\; G_z \,])}$ the corresponding smallest input containing subspaces, and by $(\hat{\gS}_i)_{i \in \mathbb{N}}$ and $(\tilde{\gS}_i)_{i \in \mathbb{N}}$ the two sequences that converge in at most $n-1$ steps to $\gS^\star_{\scriptscriptstyle (A,B,E,D_z)}$ and $\gS^\star_{\scriptscriptstyle (A,[\,B \;\; H \,],E,[\, D_z \;\; G_z \,])}$, respectively. In general, $\gV^\star_{\scriptscriptstyle (A,B,E,D_z)}\subseteq \gV^\star_{\scriptscriptstyle (A,[\,B \;\; H \,],E,[\, D_z \;\; G_z \,])}$ and $\gS^\star_{\scriptscriptstyle (A,B,E,D_z)}\subseteq \gS^\star_{\scriptscriptstyle (A,[\,B \;\; H \,],E,[\, D_z \;\; G_z \,])}$; indeed, $\hat{\gV}_i \subseteq \tilde{\gV}_i$ and $\hat{\gS}_i \subseteq \tilde{\gS}_i$ for all $i \in \mathbb{N}$. However, when the inclusion $\ima \bsmat H \\[1mm] G_z \esmat\subseteq (\gV^\star_{\scriptscriptstyle (A,B,E,D_z)}\oplus 0_{\scriptscriptstyle \gZ})+\ima  \bsmat B \\[1mm] D_z \esmat$ holds true, we have $\hat{\gV}_i = \tilde{\gV}_i$ for all $i \in \mathbb{N}$, \cite[Lemma 3]{N-TAC-08}. Even if we still have $\hat{\gS}_i \subseteq \tilde{\gS}_i$ for all $i \in \mathbb{N}$, the identity $\tilde{\gV}_i+\tilde{\gS}_j=\hat{\gV}_i+\tilde{\gS}_j=\hat{\gV}_i+\hat{\gS}_j$ holds for all $i,j \in \mathbb{N}$, as the following result shows.

\begin{lemma}
\label{lemfaticata}
Let $\ima \bsmat H \\[1mm] G_z \esmat\subseteq (\gV^\star_{\scriptscriptstyle (A,B,E,D_z)}\oplus 0_{\scriptscriptstyle \gZ})+\ima  \bsmat B \\[1mm] D_z \esmat$ hold. Then,
\[
\hat{\gV}_i+\tilde{\gS}_j=\hat{\gV}_i+\hat{\gS}_j
\]
for all $i,j\in \mathbb{N}$.
 \end{lemma}
\proof
 We start proving by induction that $\tilde{\gS}_j\subseteq \gV^\star_{\scriptscriptstyle (A,B,E,D_z)}+\hat{\gS}_j$ for all $j\in \mathbb{N}$. The statement is trivially true for $j=0$. Suppose that $\tilde{\gS}_i\subseteq \gV^\star_{\scriptscriptstyle (A,B,E,D_z)}+\hat{\gS}_i$ for a certain $i\in \mathbb{N}$, and we prove that $\tilde{\gS}_{i+1}\subseteq \gV^\star_{\scriptscriptstyle (A,B,E,D_z)}+\hat{\gS}_{i+1}$. Let $\x\in \tilde{\gS}_{i+1}$. There exist $\x_1\in \tilde{\gS}_i$, $\u\in \gU$ and $\w\in 
\gW$ such that 
$\x =A\,\x_1+B\,\u+H\,\w$ and $E\,\x_1+D_z\,\u+G_z\,\w=0$.
From $\ima \bsmat H \\[1mm] G_z \esmat\subseteq (\gV^\star_{\scriptscriptstyle (A,B,E,D_z)}\oplus 0_{\scriptscriptstyle \gZ})+\ima  \bsmat B \\[1mm] D_z \esmat$, we can find two matrices $M$ and $N$ of suitable sizes such that 
$H=V\,M+B\,N$ and $G_z=D_z\,N$, where $V$ is a basis matrix of $\gV^\star_{\scriptscriptstyle (A,B,E,D_z)}$. 
We can rewrite the previous two identities as
\[
\x = A\,\x_1+B\,\u+(V\,M+B\,N)\,\w
\]
 and
\[
E\,\x_1+D_z\,\u+(D_z\,N)\,\w=0,
\]
i.e.,
\[
\x =A\,\x_1+B\,(\u+N\,\w)+V\,M\,\w
\]
 and  
\[
E\,\x_1+D_z\,(\u+N\,\w)=0.
\]
Since $\x_1\in \tilde{\gS}_i\subseteq \gV^\star_{\scriptscriptstyle (A,B,E,D_z)}+\hat{\gS}_i$, from the inductive assumption, we can write $\x_1=\x_v+\x_s$, where $\x_v\in \gV^\star_{\scriptscriptstyle (A,B,E,D_z)}$ and $\x_s\in \hat{\gS}_i$, so that 
\[
\x= A\,\x_v+A\,\x_s+B\,(\u+N\,\w)+V\,M\,\w
\]
 and 
 \[
 E\,\x_v+E\,\x_s+D_z\,(\u+N\,\w)=0.
 \]
Let $F\in \mathfrak{F}_{\scriptscriptstyle (A,B,E,D_z)}(\gV^\star_{\scriptscriptstyle (A,B,E,D_z)})$. Adding and subtracting $B\,F\,\x_v$ in the right hand-side of the first equation and $D\,F\,\x_v$ in the right hand-side of the second equation gives
\beann
 \x &=&  A\,\x_s+B\,(\u+N\,\w-F\,\x_v)+(A+B\,F)\,\x_v+V\,M\,\w, \\
0  &=&  E\,\x_s+D_z\,(\u+N\,\w-F\,\x_v)+(E+D_z\,F)\,\x_v.
\eeann
Clearly, $(A+B\,F)\,\x_v+V\,M\,\w\in \gV^\star_{\scriptscriptstyle (A,B,E,D_z)}$ and $(E+D_z\,F)\,\x_v=0$. Defining ${\omega}=\u+N\,\w-F\,\x_v$
and ${\xi}=A\,\x_s+B\,{\omega}$, since $E\,\x_s+D_z\,{\omega}=0$ with $\x_s\in \hat{\gS}_i$, it follows that ${\xi}\in \hat{\gS}_{i+1}$. Thus, 
 $\x\in \hat{\gS}_{i+1}+\gV^\star_{\scriptscriptstyle (A,B,E,D_z)}$ as required.
 We have proved that $\tilde{\gS}_j\subseteq \gV^\star_{\scriptscriptstyle (A,B,E,D_z)}+\hat{\gS}_j$ for all $j\in \mathbb{N}$. Clearly, $\gV^\star_{\scriptscriptstyle (A,B,E,D_z)}+\tilde{\gS}_j\subseteq \gV^\star_{\scriptscriptstyle (A,B,E,D_z)}+\hat{\gS}_j$ for all $j\in \mathbb{N}$. Since $\gV^\star_{\scriptscriptstyle (A,B,E,D_z)}=\gV^\star_{\scriptscriptstyle (A,[\,B \;\; H \,],E,[\, D_z \;\; G_z \,])}$, we have 
 $\gV^\star_{\scriptscriptstyle (A,[\,B \;\; H \,],E,[\, D_z \;\; G_z \,])}+\tilde{\gS}_j\subseteq \gV^\star_{\scriptscriptstyle (A,B,E,D_z)}+\hat{\gS}_j$ for all $j\in \mathbb{N}$.
 Since we showed that $\gV^\star_{\scriptscriptstyle (A,[\,B \;\; H \,],E,[\, D_z \;\; G_z \,])}+\tilde{\gS}_j\supseteq \gV^\star_{\scriptscriptstyle (A,B,E,D_z)}+\hat{\gS}_j$, we get $\gV^\star_{\scriptscriptstyle (A,B,E,D_z)}+\tilde{\gS}_j = \gV^\star_{\scriptscriptstyle (A,B,E,D_z)}+\hat{\gS}_j$ for all $j\in \mathbb{N}$.
 Finally, since $\hat{\gV}_i \supseteq \gV^\star_{\scriptscriptstyle (A,B,E,D_z)}$ for all $i\in \mathbb{N}$, then 
$\hat{\gV}_i+\tilde{\gS}_j = \hat{\gV}_i+\hat{\gS}_j$
 for all $i,j\in \mathbb{N}$.  
\endproof
 
Following the notation of \cite{Basile-M-92}, we denote
\beann
\gV_m &\defi &\gR^\star_{\scriptscriptstyle (A,[\,B \;\; H \,],E,[\, D_z \;\; G_z \,])} \\
&=& \gV^\star_{\scriptscriptstyle (A,[\,B \;\; H \,],E,[\, D_z \;\; G_z \,])} \cap \gS^\star_{\scriptscriptstyle (A,[\,B \;\; H \,],E,[\, D_z \;\; G_z \,])}=\min \Phi_{\scriptscriptstyle (A,[\,B \;\; H \,],E,[\, D_z \;\; G_z \,])}.
\eeann
If $\ima \bsmat H \\[1mm] G_z \esmat\subseteq (\gV^\star_{\scriptscriptstyle (A,B,E,D_z)}\oplus 0_{\scriptscriptstyle \gZ})+\ima  \bsmat B \\[1mm] D_z \esmat$, then we have $\gV_m=\gV^\star_{\scriptscriptstyle (A,B,E,D_z)} \cap \gS^\star_{\scriptscriptstyle (A,[\,B \;\; H \,],E,[\, D_z \;\; G_z \,])}$ in view of Theorem~\ref{the14}.

We now consider the two quadruples $(A,H,C,G_y)$ and $\left(A,H,\bsmat C\\[0.8mm] E \esmat,\bsmat G_y \\[-0.5mm] G_z \esmat\right)$. We denote by $(\check{\gV}_i)_{i \in \mathbb{N}}$ and $(\bar{\gV}_i)_{i \in \mathbb{N}}$ the two sequences that converge in at most $n-1$ steps to $\gV^\star_{\scriptscriptstyle (A,H,C,G_y)}$ and $\gV^\star_{\scriptscriptstyle \left(A,H,\bsmat {\scriptscriptstyle C}\\[0.3mm] {\scriptscriptstyle E} \esmat,\bsmat {\scriptscriptstyle G_y} \\[-0.5mm] {\scriptscriptstyle G_z} \esmat\right)}$, respectively. Similarly, we denote by $(\check{\gS}_i)_{i \in \mathbb{N}}$ and $(\bar{\gS}_i)_{i \in \mathbb{N}}$ the two sequences that converge in at most $n-1$ steps to $\gS^\star_{\scriptscriptstyle (A,H,C,G_y)}$ and $\gS^\star_{\scriptscriptstyle \left(A,H,\bsmat {\scriptscriptstyle C}\\[0.3mm] {\scriptscriptstyle E} \esmat,\bsmat {\scriptscriptstyle G_y} \\[-0.5mm] {\scriptscriptstyle G_z} \esmat\right)}$, respectively. 
In general, $\gV^\star_{\scriptscriptstyle \left(A,H,\bsmat {\scriptscriptstyle C}\\[0.3mm] {\scriptscriptstyle E} \esmat,\bsmat {\scriptscriptstyle G_y} \\[-0.5mm] {\scriptscriptstyle G_z} \esmat\right)} \subseteq \gV^\star_{\scriptscriptstyle (A,H,C,G_y)}$ and $\gS^\star_{\scriptscriptstyle \left(A,H,\bsmat {\scriptscriptstyle C}\\[0.3mm] {\scriptscriptstyle E} \esmat,\bsmat {\scriptscriptstyle G_y} \\[-0.5mm] {\scriptscriptstyle G_z} \esmat\right)} \subseteq \gS^\star_{\scriptscriptstyle (A,H,C,G_y)}$.
When the inclusion $\ker [\begin{array}{cc} E & G_z \end{array}]\supseteq (\gS^\star\oplus \gW)\cap \ker [\begin{array}{cc} C & G_y \end{array}]$ holds, we have $\check{\gS}_i = \bar{\gS}_i$ for all $i \in \mathbb{N}$ from the dual of \cite[Lemma 3]{N-TAC-08}. 
The following result can be proved by dualizing the proof of Lemma~\ref{lemfaticata}.

\begin{lemma}
Let $\ker [\begin{array}{cc} E & G_z \end{array}]\supseteq (\gS^\star_{\scriptscriptstyle (A,H,C,G_y)} \oplus \gW)\cap \ker [\begin{array}{cc} C & G_y \end{array}]$. For all $i,j\in \mathbb{N}$ there holds
\[
\bar{\gV}_i\cap \check{\gS}_j=\check{\gV}_i\cap \check{\gS}_j.
\]
 \end{lemma}

Following the notation of \cite{Basile-M-92}, we denote
\beann
\gS_M &\defi &\gQ^\star_{\scriptscriptstyle \left(A,H,\bsmat {\scriptscriptstyle C}\\[0.3mm] {\scriptscriptstyle E} \esmat,\bsmat {\scriptscriptstyle G_y} \\[-0.5mm] {\scriptscriptstyle G_z} \esmat\right)}\\
&=& \gV^\star_{\scriptscriptstyle \left(A,H,\bsmat {\scriptscriptstyle C}\\[0.3mm] {\scriptscriptstyle E} \esmat,\bsmat {\scriptscriptstyle G_y} \\[-0.5mm] {\scriptscriptstyle G_z} \esmat\right)} + \gS^\star_{\scriptscriptstyle \left(A,H,\bsmat {\scriptscriptstyle C}\\[0.3mm] {\scriptscriptstyle E} \esmat,\bsmat {\scriptscriptstyle G_y} \\[-0.5mm] {\scriptscriptstyle G_z} \esmat\right)}=\max \Psi_{\scriptscriptstyle \left(A,H,\bsmat {\scriptscriptstyle C}\\[0.3mm] {\scriptscriptstyle E} \esmat,\bsmat {\scriptscriptstyle G_y} \\[-0.5mm] {\scriptscriptstyle G_z} \esmat\right)}.
\eeann
If $\ker [\begin{array}{cc} E & G_z \end{array}]\supseteq (\gS^\star_{\scriptscriptstyle (A,H,C,G_y)} \oplus \gW)\cap \ker [\begin{array}{cc} C & G_y \end{array}]$, we have $\gS_M=\gV^\star_{\scriptscriptstyle \left(A,H,\bsmat {\scriptscriptstyle C}\\[0.3mm] {\scriptscriptstyle E} \esmat,\bsmat {\scriptscriptstyle G_y} \\[-0.5mm] {\scriptscriptstyle G_z} \esmat\right)} + \gS^\star_{\scriptscriptstyle (A,H,C,G_y)}$.

The proof of the following result is straightforward.

\begin{lemma}
\label{gt}
The following inclusions hold:
\begin{itemize}
\item $\gV^\star_{\scriptscriptstyle \left(A,H,\bsmat {\scriptscriptstyle C}\\[0.3mm] {\scriptscriptstyle E} \esmat,\bsmat {\scriptscriptstyle G_y} \\[-0.5mm] {\scriptscriptstyle G_z} \esmat\right)}\subseteq  \gV^\star_{\scriptscriptstyle (A,B,E,D_z)}\subseteq 
\gV^\star_{\scriptscriptstyle (A,[\,B \; H \,],E,[\, D_z \; G_z \,])}$;
\item $\gS^\star_{\scriptscriptstyle \left(A,H,\bsmat {\scriptscriptstyle C}\\[0.3mm] {\scriptscriptstyle E} \esmat,\bsmat {\scriptscriptstyle G_y} \\[-0.5mm] {\scriptscriptstyle G_z} \esmat\right)}\subseteq \gS^\star_{\scriptscriptstyle (A,H,C,G_y)} \subseteq \gS^\star_{\scriptscriptstyle (A,[\,B \; H \,],E,[\, D_z \; G_z \,])}$;
\item $\gS^\star_{\scriptscriptstyle (A,H,C,G_y)} \subseteq \gV^\star_{\scriptscriptstyle (A,B,E,D_z)} \;\Rightarrow \;\gS^\star_{\scriptscriptstyle \left(A,H,\bsmat {\scriptscriptstyle C}\\[0.3mm] {\scriptscriptstyle E} \esmat,\bsmat {\scriptscriptstyle G_y} \\[-0.5mm] {\scriptscriptstyle G_z} \esmat\right)} \subseteq \gV^\star_{\scriptscriptstyle (A,[\,B \; H \,],E,[\, D_z \; G_z \,])}$.
\end{itemize}
\end{lemma}

\begin{lemma}
\label{lem22}
Let $ \gS^\star_{\scriptscriptstyle (A,H,C,G_y)} \subseteq \gV^\star_{\scriptscriptstyle (A,B,E,D_z)}$. Then,
the subspace $\gV_m+\gS_M$ is $(A,[\begin{array}{cc} B & H \end{array}],E,[\begin{array}{cc} D_z & G_z\end{array}])$-self bounded, and 
the subspace $\gV_m\cap \gS_M$ is $\left(A,H ,\bsmat C\\[1mm] E\esmat,\bsmat G_y \\[1mm] G_z\esmat\right)$-self hidden. 
\end{lemma}

\proof
We find
\bea
 \gV_m&+&\gS_M = (\gV^\star_{\scriptscriptstyle (A,[\,B \; H \,],E,[\, D_z \; G_z \,])}\cap \gS^\star_{\scriptscriptstyle (A,[\,B \; H \,],E,[\, D_z \; G_z \,])}) \nonumber \\
&& +(\gS^\star_{\scriptscriptstyle \left(A,H,\bsmat {\scriptscriptstyle C}\\[0.3mm] {\scriptscriptstyle E} \esmat,\bsmat {\scriptscriptstyle G_y} \\[-0.5mm] {\scriptscriptstyle G_z} \esmat\right)}+\gV^\star_{\scriptscriptstyle \left(A,H,\bsmat {\scriptscriptstyle C}\\[0.3mm] {\scriptscriptstyle E} \esmat,\bsmat {\scriptscriptstyle G_y} \\[-0.5mm] {\scriptscriptstyle G_z} \esmat\right)}) \nonumber \\
& = &  \bigl((\gV^\star_{\scriptscriptstyle (A,[\,B \; H \,],E,[\, D_z \; G_z \,])}\cap \gS^\star_{\scriptscriptstyle (A,[\,B \;H \,],E,[\, D_z \; G_z \,])})+\gS^\star_{\scriptscriptstyle \left(A,H,\bsmat {\scriptscriptstyle C}\\[0.3mm] {\scriptscriptstyle E} \esmat,\bsmat {\scriptscriptstyle G_y} \\[-0.5mm] {\scriptscriptstyle G_z} \esmat\right)}\bigr)  \nonumber \\
&& +\gV^\star_{\scriptscriptstyle \left(A,H,\bsmat {\scriptscriptstyle C}\\[0.3mm] {\scriptscriptstyle E} \esmat,\bsmat {\scriptscriptstyle G_y} \\[-0.5mm] {\scriptscriptstyle G_z} \esmat\right)}  \nonumber \\
& = &  \bigl[(\gV^\star_{\scriptscriptstyle (A,[\,B \; H \,],E,[\, D_z \; G_z \,])}+ \gS^\star_{\scriptscriptstyle \left(A,H,\bsmat {\scriptscriptstyle C}\\[0.3mm] {\scriptscriptstyle E} \esmat,\bsmat {\scriptscriptstyle G_y} \\[-0.5mm] {\scriptscriptstyle G_z} \esmat\right)})\cap (\gS^\star_{\scriptscriptstyle (A,[\,B \; H \,],E,[\, D_z \; G_z \,])} \nonumber \\
&&
+ \gS^\star_{\scriptscriptstyle \left(A,H,\bsmat {\scriptscriptstyle C}\\[0.3mm] {\scriptscriptstyle E} \esmat,\bsmat {\scriptscriptstyle G_y} \\[-0.5mm] {\scriptscriptstyle G_z} \esmat\right)})\bigr] +\gV^\star_{\scriptscriptstyle \left(A,H,\bsmat {\scriptscriptstyle C}\\[0.3mm] {\scriptscriptstyle E} \esmat,\bsmat {\scriptscriptstyle G_y} \\[-0.5mm] {\scriptscriptstyle G_z} \esmat\right)}
\nonumber \\
& = &  (\gV^\star_{\scriptscriptstyle (A,[\,B \; H \,],E,[\, D_z \; G_z \,])}\cap \gS^\star_{\scriptscriptstyle (A,[\,B \; H \,],E,[\, D_z \; G_z \,])})+\gV^\star_{\scriptscriptstyle \left(A,H,\bsmat {\scriptscriptstyle C}\\[0.3mm] {\scriptscriptstyle E} \esmat,\bsmat {\scriptscriptstyle G_y} \\[-0.5mm] {\scriptscriptstyle G_z} \esmat\right)} 
  \nonumber  \\
& = &  \gV_m+\gV^\star_{\scriptscriptstyle \left(A,H,\bsmat {\scriptscriptstyle C}\\[0.3mm] {\scriptscriptstyle E} \esmat,\bsmat {\scriptscriptstyle G_y} \\[-0.5mm] {\scriptscriptstyle G_z} \esmat\right)}, \label{eq:pippo}
\eea
in view of the modular rule \cite[p.~16]{Trentelman-SH-01} and  Lemma~\ref{gt}.
We show that $\gV_m+\gS_M$ is $(A,[\begin{array}{cc} B & H \end{array}],E,[\begin{array}{cc} D_z & G_z\end{array}])$-output nulling. The inclusion
\[
\bmat{c}   A  \\[0mm]  C  \\[0mm]  E  \emat \gV^\star_{\scriptscriptstyle \left(A,H,\bsmat {\scriptscriptstyle C}\\[0.3mm] {\scriptscriptstyle E} \esmat,\bsmat {\scriptscriptstyle G_y} \\[-0.5mm] {\scriptscriptstyle G_z} \esmat\right)}    \subseteq (\gV^\star_{\scriptscriptstyle \left(A,H,\bsmat {\scriptscriptstyle C}\\[0.3mm] {\scriptscriptstyle E} \esmat,\bsmat {\scriptscriptstyle G_y} \\[-0.5mm] {\scriptscriptstyle G_z} \esmat\right)}   \oplus 0_{\scriptscriptstyle  \gY \oplus \gZ})+\ima \bmat{c}   H   \\[0mm]   G_y   \\[0mm]   G_z   \emat
\]
implies $\bsmat  A \\[1mm]   E  \esmat  \gV^\star_{\scriptscriptstyle \left(A,H,\bsmat {\scriptscriptstyle C}\\[0.3mm] {\scriptscriptstyle E} \esmat,\bsmat {\scriptscriptstyle G_y} \\[-0.5mm] {\scriptscriptstyle G_z} \esmat\right)}    \subseteq   (\gV^\star_{\scriptscriptstyle \left(A,H,\bsmat {\scriptscriptstyle C}\\[0.3mm] {\scriptscriptstyle E} \esmat,\bsmat {\scriptscriptstyle G_y} \\[-0.5mm] {\scriptscriptstyle G_z} \esmat\right)}   \oplus 0_{\scriptscriptstyle \gZ})+\ima \bsmat H \\[1mm] G_z \esmat$, which in turn leads to
$\bsmat A\\[1mm] E\esmat\,\gV^\star_{\scriptscriptstyle \left(A,H,\bsmat {\scriptscriptstyle C}\\[0.3mm] {\scriptscriptstyle E} \esmat,\bsmat {\scriptscriptstyle G_y} \\[-0.5mm] {\scriptscriptstyle G_z} \esmat\right)}  \subseteq (\gV^\star_{\scriptscriptstyle \left(A,H,\bsmat {\scriptscriptstyle C}\\[0.3mm] {\scriptscriptstyle E} \esmat,\bsmat {\scriptscriptstyle G_y} \\[-0.5mm] {\scriptscriptstyle G_z} \esmat\right)}  \oplus 0_{\scriptscriptstyle \gZ})+\ima \bsmat B && H \\[1mm] D_z && G_z \esmat$.
Adding this to $\bsmat A\\[1mm] E\esmat\,\gV_m \subseteq (\gV_m\oplus 0_{\scriptscriptstyle \gZ})+\ima \bsmat B && H \\[1mm] D_z && G_z \esmat$ (since $\gV_m\in \Phi_{\scriptscriptstyle (A,[\, B \;\;H \,],E,[\, D_z \;\; G_z \,]}$) yields
\bea
&&\bmat{c} A\\[-1mm] E\emat\,(\gV_m+\gV^\star_{\scriptscriptstyle \left(A,H,\bsmat {\scriptscriptstyle C}\\[0.3mm] {\scriptscriptstyle E} \esmat,\bsmat {\scriptscriptstyle G_y} \\[-0.5mm] {\scriptscriptstyle G_z} \esmat\right)} )\nonumber  \\ && \qquad \subseteq \bigl((\gV_m+\gV^\star_{\scriptscriptstyle \left(A,H,\bsmat {\scriptscriptstyle C}\\[0.3mm] {\scriptscriptstyle E} \esmat,\bsmat {\scriptscriptstyle G_y} \\[-0.5mm] {\scriptscriptstyle G_z} \esmat\right)} )\oplus 0_{\scriptscriptstyle \gZ}\bigr) +\ima \bmat{cc} B & H \\[-1mm] D_z & G_z \emat.\label{last}
\eea
Thus, $\gV_m+\gS_M$ is $(A,[\begin{array}{cc} B & H \end{array}],E,[\begin{array}{cc} D_z & G_z\end{array}])$-output nulling. 
The fact that $\gV_m+\gV^\star_{\scriptscriptstyle \left(A,H,\bsmat {\scriptscriptstyle C}\\[0.3mm] {\scriptscriptstyle E} \esmat,\bsmat {\scriptscriptstyle G_y} \\[-0.5mm] {\scriptscriptstyle G_z} \esmat\right)}$ is self bounded follows immediately from the inclusion
$\gV_m+\gV^\star_{\scriptscriptstyle \left(A,H,\bsmat {\scriptscriptstyle C}\\[0.3mm] {\scriptscriptstyle E} \esmat,\bsmat {\scriptscriptstyle G_y} \\[-0.5mm] {\scriptscriptstyle G_z} \esmat\right)} \supseteq\gV_m \supseteq  \gV^\star_{\scriptscriptstyle (A,[\,B\;\;H\,],E,[\,D_z\;\; G_z\,])} \cap [\begin{array}{cc} B & H \end{array}]\,\ker [\begin{array}{cc} D_z & G_z \end{array}]$.
The second statement follows by duality. \endproof

\begin{corollary}
\label{cor5}
Let $ \gS^\star_{\scriptscriptstyle (A,H,C,G_y)} \subseteq \gV^\star_{\scriptscriptstyle (A,B,E,D_z)}$. 
The following results hold:\\
$\bullet \quad$  If $\ima \bsmat H \\[1mm] G_z \esmat \subseteq (\gV^\star_{\scriptscriptstyle (A,B,E,D_z)} \oplus 0_{\scriptscriptstyle \gZ})+\ima \bsmat B\\[1mm] D_z \esmat$, 
then $\gV_m+\gS_M$  is $(A,B,E,D_z)$-self bounded. \\
$\bullet \quad$  If  $\ker \,[\begin{array}{cc} E & G_z \end{array}] \supseteq (\gS^\star_{\scriptscriptstyle (A,H,C,G_y)} \oplus \gW)  \cap \ker \,[\begin{array}{cc}  C & G_y \end{array}]$, then $\gV_m\cap \gS_M$ is $(A,H,C,G_y)$-self hidden.
\end{corollary}
\proof
Recall that $\ima \bsmat H \\[1mm] G_z \esmat \subseteq (\gV^\star_{\scriptscriptstyle (A,B,E,D_z)} \oplus 0_{\scriptscriptstyle \gZ})+\ima \bsmat B\\[1mm] D_z \esmat$ implies $\ima \bsmat H \\[1mm] G_z \esmat \subseteq (\gV_m \oplus 0_{\scriptscriptstyle \gZ})+\ima \bsmat B\\[1mm] D_z \esmat$ from Theorem~\ref{the14}.
Using this inclusion into (\ref{last}) we obtain
\beann
&& \bmat{c} A\\[-1mm] E\emat\,(\gV_m+\gV^\star_{\scriptscriptstyle \left(A,H,\bsmat {\scriptscriptstyle C}\\[0.3mm] {\scriptscriptstyle E} \esmat,\bsmat {\scriptscriptstyle G_y} \\[-0.5mm] {\scriptscriptstyle G_z} \esmat\right)}) \\
&&\;\; \subseteq  \bigl((\gV_m+\gV^\star_{\scriptscriptstyle \left(A,H,\bsmat {\scriptscriptstyle C}\\[0.3mm] {\scriptscriptstyle E} \esmat,\bsmat {\scriptscriptstyle G_y} \\[-0.5mm] {\scriptscriptstyle G_z} \esmat\right)})\oplus 0_{\scriptscriptstyle \gZ}\bigr) +\ima \bmat{cc} B & H \\[-1mm]  D_z & G_z \emat \\
&&\;\; =\bigl((\gV_m+\gV^\star_{\scriptscriptstyle \left(A,H,\bsmat {\scriptscriptstyle C}\\[0.3mm] {\scriptscriptstyle E} \esmat,\bsmat {\scriptscriptstyle G_y} \\[-0.5mm] {\scriptscriptstyle G_z} \esmat\right)})\oplus 0_{\scriptscriptstyle \gZ}\bigr) +\ima \bmat{cc} B \\[-1mm]  D_z \emat+\ima \bmat{c} H \\[-1mm]  G_z \emat \\
&&\;\; \subseteq \bigl((\gV_m+\gV^\star_{\scriptscriptstyle \left(A,H,\bsmat {\scriptscriptstyle C}\\[0.3mm] {\scriptscriptstyle E} \esmat,\bsmat {\scriptscriptstyle G_y} \\[-0.5mm] {\scriptscriptstyle G_z} \esmat\right)})\oplus 0_{\scriptscriptstyle \gZ}\bigr) +\ima \bmat{cc} B \\[-1mm]  D_z \emat+(\gV_m \oplus 0_{\scriptscriptstyle \gZ})\\ 
&& \;\;= \bigl((\gV_m+\gV^\star_{\scriptscriptstyle \left(A,H,\bsmat {\scriptscriptstyle C}\\[0.3mm] {\scriptscriptstyle E} \esmat,\bsmat {\scriptscriptstyle G_y} \\[-0.5mm] {\scriptscriptstyle G_z} \esmat\right)})\oplus 0_{\scriptscriptstyle \gZ}\bigr) +\ima \bmat{cc} B \\[-1mm]  D_z \emat.
\eeann
We also need to prove that $\gV_m+\gS_M\supseteq \gV^\star_{\scriptscriptstyle (A,B,E,D_z)}  \cap B\,\ker D_z$: this follows from
$\gV_m+\gS_M\supseteq \gV^\star_{\scriptscriptstyle (A,B,E,D_z)} \cap [\begin{array}{cc} B & H \end{array}]\,\ker [\begin{array}{cc} D_z & G_z \end{array}]$.
The second can be proved by duality.
\endproof

\section{Problem solution}
We begin by first presenting the following result, see \cite[Lemma 3.2]{Stoorvogel-W-91}. The proof can be carried out along the same lines of the proof of \cite[Lemma 6.3]{Trentelman-SH-01}. The next few preliminary results involve integers $n_1,n_2,m,p\in \mathbb{N} \setminus \{0\}$, a field $\field$, a subspace
$\gM$ of $\field^{n_2}$ and a subspace $\gN$ of $\field^{n_1}$.
We also consider the matrices $\tilde{A}\in \field^{n_1 \times n_2}$, $\tilde{B}\in \field^{n_1 \times m}$ and $\tilde{C}\in \field^{p\times n_2}$.

\begin{lemma}
\label{lemstoor}
There holds
$\tilde{A}\,\gM \subseteq \gN+\ima \tilde{B}$ and $\tilde{A}\,(\gM\cap \ker \tilde{C}) \subseteq \gN$
if and only if there exists  $K\in \real^{m \times p}$ such that
$(\tilde{A}+\tilde{B}\,K\,\tilde{C})\,\gM\subseteq \gN$.
\end{lemma}

\begin{lemma}
\label{lem3cond}
Let $\gV$ be an $(A,B,E,D_z)$-output nulling subspace and let $\gS$ be an $(A,H,C,G_y)$-input containing subspace. If\\[-.7cm]
\begin{description}
\item{\bf (a)} $\ima \bmat{c} H \\[-1.5mm] G_z \emat\subseteq (\gV \oplus 0_{\scriptscriptstyle \gZ}) +\ima \bmat{c} B \\[-1.5mm] D_z\emat$,
\item{\bf (b)} $\ker \,[\begin{array}{cc} E & G_z \end{array}] \supseteq (\gS \oplus \gW)  \cap \ker \,[\begin{array}{cc}  C & G_y \end{array}] $,
\item{\bf (c)} $\gS\subseteq \gV$,
\end{description}
then there exists an output feedback matrix $K$ such that 
 \bea
 \label{clK}
 \bmat{cc} A+B\,K\,C & H+B\,K\,G_y \\[-1.5mm]
 E+D_z\,K\,C & G_z+ D_z\,K\,G_y \emat(\gS \oplus \gW)\subseteq \gV \oplus 0_{\scriptscriptstyle \gZ}.
 \eea
 Conversely, if $K$ exists such that (\ref{clK}) holds, then {\bf (a-b)} hold. 
\end{lemma}
\proof
We prove that if {\bf (a-c)} hold, then $K$ exists such that (\ref{clK}) holds.
 Since $\gV$ is $(A,B,E,D_z)$-output nulling, we have
 $\bsmat A \\[1mm] E \esmat\,\gV\subseteq \gV \oplus 0_{\scriptscriptstyle \gZ} +\ima \bsmat B \\[1mm] D_z\esmat$.
  Combining this inclusion with {\bf (a)} yields
 $\bsmat A && H \\[1mm] E && G_z \esmat\,(\gV \oplus \gW) \subseteq \gV \oplus 0_{\scriptscriptstyle \gZ} +\ima \bsmat B \\[1mm] D_z\esmat$.
From {\bf (c)}, we also have
  \bea
  \label{inizio}
 \bmat{cc} A & H \\[-1.5mm] E & G_z \emat\,(\gS\oplus \gW) \subseteq \gV \oplus 0_{\scriptscriptstyle \gZ} +\ima \bmat{c} B \\[-1.5mm] D_z\emat.
 \eea
 Similarly, since $\gS$ is $(A,H,C,G_y)$-input containing, we have
 $[\begin{array}{cc} A & H \end{array}] \left(\gS \oplus \gW \cap \ker  [\begin{array}{cc}  C & G_y \end{array}] \right) \subseteq \gS$.
Taking {\bf (b)} into account gives
 $\bsmat A && H \\[1mm] E && G_z \esmat\,\left(\gS \oplus \gW \cap \ker  [\begin{array}{cc}  C & G_y \end{array}] \right)\subseteq 
 \gS \oplus 0_{\scriptscriptstyle \gZ}$. Again, from {\bf (c)} we obtain
 \bea
 \label{fine}
 \bmat{cc} A & H \\[-1.5mm] E & G_z \emat\,\left(\gS \oplus \gW \cap \ker  [\begin{array}{cc}  C & G_y \end{array}] \right)\subseteq 
  \gV \oplus 0_{\scriptscriptstyle \gZ}.
 \eea
We can now apply Lemma~\ref{lemstoor} considering the two inclusions (\ref{inizio}) and (\ref{fine}), i.e., by considering
$\tilde{A}\rightarrow\bsmat A && H \\[1mm] E && G_z \esmat$, $\tilde{B} \rightarrow\bsmat B \\[1mm] D_z\esmat$,  $\tilde{C} \rightarrow  [\begin{array}{cc}  C & G_y \end{array}]$, as well as the subspaces $\gM=\gS \oplus \gW$ and $\gN=\gV \oplus 0_{\scriptscriptstyle \gZ}$.  Thus,
 there exists $K\in \real^{p \times m}$ such that
 \[
 \left(\bmat{cc} A & H \\[-1.5mm] E & G_z \emat+\bmat{c} B \\[-1.5mm] D_z\emat\,K\, [\begin{array}{cc} C & G_y \end{array}] \right)(\gS \oplus \gW)\subseteq \gV \oplus 0_{\scriptscriptstyle \gZ},
 \]
which is exactly (\ref{clK}).
We now prove the converse.
 Let $K$ be such that (\ref{clK}) holds.
 Let $S$ be a basis matrix of $\gS$ and $V$ be a basis matrix of $\gV$. We can re-write (\ref{clK}) as
 \bea
 \label{h2}
 \bmat{cc} A+B\,K\,C & H+B\,K\,G_y \\[-1.5mm]
 E+D_z\,K\,C & G_z+ D_z\,K\,G_y \emat\bmat{cc} S & 0 \\[-1.5mm] 0 & I \emat=\bmat{c} V \\[-1.5mm] 0 \emat\,X
 \eea
for some matrix $X$ of suitable size, which gives the two equations
$H+B\,K\,G_y=V\,X$ and
$G_z+D_z\,K\,G_y=0$. These can be rewritten together as
$\bsmat H \\[1mm] G_z \esmat=\bsmat V \\[1mm] 0 \esmat\,X+\bsmat B \\[1mm] D_z \esmat\,(-K\,G_y)$,
 so that {\bf (a)} holds.
 From (\ref{h2}) we also find
 \bea
 \label{poi}
[\begin{array}{cc} E & G_z\end{array}]\bmat{cc} S & 0 \\[-1.5mm] 0 & I \emat+D_z\,K\,[\begin{array}{cc} C & D_y\end{array}]\bmat{cc} S & 0 \\[-1.5mm] 0 & I \emat=0.
 \eea
 Let $\bsmat \x \\[1mm] \y\esmat\in \gS\oplus \gW\cap \ker [\begin{array}{cc} C & D_y\end{array}]$. Then there exists $\eta$ such that $\bsmat \x \\[1mm] \y\esmat=\bsmat S & 0 \\[1mm] 0 & I \esmat\,\eta$. Multiplying (\ref{poi}) by $\eta$, since $\bsmat \x \\[1mm] \y\esmat\in \ker [\begin{array}{cc} C & D_y\end{array}]$, we find
 $[\begin{array}{cc} E & G_z\end{array}]\bsmat \x \\[1mm] \y \esmat=0$,
 so that {\bf (b)} holds.
\endproof

\begin{example}
{
The existence of a matrix $K$ satisfying (\ref{clK}) for an $(A,B,E,D_z)$-output nulling subspace $\gV$ and an $(A,H,C,G_y)$-input containing subspace $\gS$ does not imply the condition $\gS\subseteq \gV$. Consider for example
\beann
A \ns&\ns = \ns&\ns\bsmat 0 &0 & 0 \\[1mm]
0 & 0 & 0 \\[1mm]
0 & -1 & 1\esmat,\;\;B=\bsmat -1 && 0\\[1mm] 1 && 0 \\[1mm] 0 && 0  \esmat,\;\;H=\bsmat 1  \\[1mm] -1  \\[1mm] 0  \esmat, \\
 C \ns&\ns = \ns&\ns \bsmat 1 & 1 & 0 \\[1mm] 0 & 1 & 0 \esmat,\;\;D_y=\bsmat 0 && 0 \\[1mm] 0 && 0 \esmat,\;\; G_y=\bsmat -1\\[1mm] 0 \esmat,\\
  E\ns&\ns = \ns&\ns \bsmat 1& 0 & 0 \esmat,\;\; D_z=[\,0\;\;1\,],\;\; G_z=1,
\eeann
with the subspaces $\gV=\spanR\left\{\bsmat 1\\0\\0\esmat,\bsmat 0\\0\\1\esmat\right\}$ and $\gS=\spanR\left\{\bsmat 0\\2\\7\esmat\right\}$. One can easily verify that $K=\bsmat 1 && -1 \\[1mm] 1 && -1 \esmat$ satisfies  (\ref{clK}), and that $\gV$ and $\gS$ are, respectively, $(A,B,E,D_z)$-output nulling  and $(A,H,C,G_y)$-input containing; in addition, $\gV$ satisfies {\bf (a)} and $\gS$ satisfies {\bf (b)} of Lemma \ref{lem3cond}. However, clearly {\bf (c)} is not satisfied in this case.
}
\end{example}

The following result contains the generalization of a fundamental property to the case where all the feedthrough matrices are allowed to be nonzero. The major technical difficulty is the fact that in this case, the well-posedness needs to be taken into account. In other words, while showing that the conditions of the following theorem are sufficient for the existence of a decoupling filter only requires more convoluted matrix manipulations with respect to the strictly proper case, the necessity needs to be addressed more carefully.

\begin{theorem}
\label{solv}
Problem \ref{pro1} is solvable if and only if there exist an $(A,B,E,D_z)$-output nulling subspace $\gV$, an $(A,H,C,G_y)$-input containing subspace $\gS$ and a matrix $K\in \real^{m \times p}$ such that
\begin{description}
\item{\bf (i)} $\ima \bmat{c} H \\[-1.5mm] G_z \emat \subseteq (\gV \oplus 0_{\scriptscriptstyle \gZ})+\ima \bmat{c} B \\[-1.5mm] D_z \emat$;
\item{\bf (ii)} $\ker \,[\begin{array}{cc} E & G_z \end{array}] \supseteq (\gS \oplus \gW)  \cap \ker \,[\begin{array}{cc}  C & G_y \end{array}] $;
\item{\bf (iii)} $\gS\subseteq \gV$;
\item{\bf (iv)} $I+K\,D_y$ is non-singular, and $K$ satisfies 
\bea
\label{clK1}
 \bmat{cc} A+B\,K\,C & H+B\,K\,G_y \\[-1.5mm]
 E+D_z\,K\,C & G_z+ D_z\,K\,G_y \emat(\gS \oplus \gW)\subseteq \gV \oplus 0_{\scriptscriptstyle \gZ}.
\eea
\end{description}
\end{theorem}

\proof
{\bf (If)}.  We define the compensator matrices as 
 \bea
 \label{eqcomp}
 \begin{array}{rclccrcl}
 A_c \ns&\ns = \ns&\ns A+G\,C+(B+G\,D_y)(I+K\,D_y)^{-1}(F-K\,C),\\
B_c \ns&\ns = \ns&\ns (B+G\,D_y)(I+K\,D_y)^{-1} K-G, \\
C_c \ns&\ns = \ns&\ns (I+K\,D_y)^{-1} (F-K\,C),\\
D_c \ns&\ns = \ns&\ns (I+K\,D_y)^{-1} \,K.\end{array}
\eea
where $F\in \mathfrak{F}_{\scriptscriptstyle (A,B,E,D_z)}(\gV)$, so that $\bsmat A+B\,F \\[1mm] E+D_z\,F\esmat \gV \subseteq \gV\oplus 0_{\scriptscriptstyle \gZ}$, and where $G\in \mathfrak{G}_{\scriptscriptstyle (A,H,C,G_y)}(\gS)$, so that $[\begin{array}{cc} A+G\,C & H+G\,G_y \end{array}] (\gS\oplus \gW) \subseteq \gS$. Using these matrices in (\ref{cl}) and using the matrix inversion lemma\footnote{Given matrices $P,Q,R,S$ of conformable sizes such that $P$, $R$ and $P+Q\,R\,S$ are invertible, there holds
$(P+Q\,R\,S)^{-1}=P^{-1}-P^{-1}\,Q\,(R^{-1}+S\,P^{-1}\,Q)^{-1}\,S\,P^{-1}$.
}, after some lengthy but standard matrix manipulations we obtain
{\beann
\widehat{A}
\ns&\ns = \ns&\ns \bmat{cc}   A+B\,K\,C  & B\,(F-K\,C)   \\[-1.5mm]
  (B\,K-G)\,C & A+G\,C+B\,F-B\,K\,C  \emat    ,\quad       \widehat{H}=
\bmat{c}   H+B\,K\,G_y   \\[-1.5mm]
  (B\,K-G)\,G_y    \emat    , \\
\widehat{C}
\ns&\ns = \ns&\ns[\begin{array}{cc}
E+D_z\,K\,C  & D_z\,(F-K\,C) \end{array}],\quad \widehat{G}=G_z+D_z\,K\,G_y.
\eeann}
Defining $\e=\x-\p$, we obtain
 \beann
\bmat{c}  \gD {\x}(t) \\[-1.5mm] \gD {\e}(t)\emat 
&  =  & \bmat{cc}  A+B\,F & B\,(K\,C-F) \\[-1.5mm]
0 & A+G\,C  \emat      \bmat{c} \x(t) \\[-1.5mm]  \e(t) \emat+
\bmat{c} H+B\,K\,G_y \\[-1.5mm]
H+G\,G_y  \emat     \w(t), \\
\z(t)
&  =  &[\begin{array}{cc}
E+D_z F  & D_z (K C-F) \end{array}] \bmat{c} \x(t)\\[-1.5mm]  \e(t) \emat+\bigl(G_z+D_z K G_y\bigr) \w(t).
\eeann
We now show that the transfer function $G_{\z,\w}(\l)$ is zero: 
{\small
\beann
&& \hspace{-1.2cm} G_{\z,\w}(\l) =  [\begin{array}{cc}
E+D_z F  & D_z (K C-F) \end{array}]  \bmat{cc}  \l I-A-B F & -B (K C-F)  \\[-1.5mm]
 0 &\l I-A-G C  \emat^{-1} \bmat{c}  H+B K G_y \\[-1.5mm]
H+G G_y  \emat+G_z+D_z\,K\,G_y \\
 &&= 
(E+D_z F)(\l I-A-B F)^{\text{\tiny $-1$}} (H+B K G_y)\\
\ns&\ns \ns&\ns 
+(E+D_z F)((\l I-A-B F))^{\text{\tiny $-1$}}(B K C-B F)(\l I-A-G C)^{\text{\tiny $-1$}}(H+G G_y)\\
\ns&\ns \ns&\ns 
 +D_z (K C-F)(\l I-A-G C)^{\text{\tiny $-1$}}(H+G G_y)+G_z+D_z K G_y \\
&&=  (E+D_z F)(\l I-A-B F)^{\text{\tiny $-1$}} (H+B K G_y) \\
\ns&\ns \ns&\ns 
+(E+D_z F){(\l I-A-B F)^{\text{\tiny $-1$}}}{(\l I-A-B F)}(\l I-A-G C)^{\text{\tiny $-1$}}(H+G G_y) \\
\ns&\ns \ns&\ns 
-(E+D_z F)(\l I-A-B F)^{\text{\tiny $-1$}}(\l I-A-B K C)(\l I-A-G C)^{\text{\tiny $-1$}}(H+G G_y)
\\
\ns&\ns \ns&\ns 
 +D_z (K C-F)(\l I-A-G C)^{\text{\tiny $-1$}}(H+G G_y)+G_z+D_z K G_y \\
&&= (E+D_z F)(\l I-A-B F)^{\text{\tiny $-1$}} (H+B K G_y)+
E (\l I-A-G C)^{\text{\tiny $-1$}} (H+G G_y)\\
\ns&\ns \ns&\ns 
-(E+D_z F)(\l I-A-B F)^{\text{\tiny $-1$}} (\l I-A-B K C)(\l I-A-G C)^{\text{\tiny $-1$}} (H+G G_y)
\\
\ns&\ns \ns&\ns 
+D_z K C(\l I-A-G C)^{\text{\tiny $-1$}}(H+G G_y)+G_z+D_z K G_y \\
&&= \underbrace{(E+D_z F)(\l I-A-B F)^{\text{\tiny $-1$}} (H+B K G_y)}_{T_1(\l)}+
\underbrace{(E+D_z K C) (\l I-A-G C)^{\text{\tiny $-1$}} (H+G G_y)}_{T_2(\l)}\\
\ns&\ns \ns&\ns 
-\underbrace{(E+D_z F)(\l I-A-B F)^{\text{\tiny $-1$}}(\l I-A-B K C)(\l I-A-G C)^{\text{\tiny $-1$}} (H+G G_y)}_{T_3(\l)}
+\underbrace{G_z+D_z K G_y}_{T_4}, 
 \eeann
 }
where we have used the identity $B\,K\,C-B\,F=(\l\,I-A-B\,F)-(\l\,I-A-B\,K\,C)$. 
Now, (\ref{clK1}) is equivalent to  
\bea
(A+B\,K\,C)\,\gS \subseteq \gV, \label{eqa0}\\
(E+D_z\,K\,C)\,\gS=0_{\scriptscriptstyle \gZ}, \label{eqb0}
\eea
\bea
\ima (H+B\,K\,G_y)\subseteq \gV, \label{eqc0}\\
G_z+D_z\,K\,G_y =0. \label{eqd0}
 \eea
Eq. (\ref{eqc0}), together with the inclusion 
\bea
\label{u1}
\ker\bigl((E+D_z\,F)\,(\l\,I-A-B\,F)^{-1}\bigr)\supseteq \gV,
\eea
 see (\ref{eqV}), yields
$\ima(E+D_z\,F)(\l\,I-A-B\,F)^{-1} (H+B\,K\,G_y) \subseteq  (E+D_z\,F)(\l\,I-A-B\,F)^{-1}\,\gV=0_{\scriptscriptstyle \gZ}$, which proves that $T_1(\l)$ is zero. Similarly, (\ref{eqb0}) with 
 \bea
 \label{u2}
 \ima \bigl((\l\,I-A-G\,C)^{-1}(H+G\,G_y)\bigr)\subseteq \gS,
 \eea
 see (\ref{eqS}),
  yields $(E+D_z\,K\,C)\,(\l\,I-A-G\,C)^{-1} (H+G\,G_y)\subseteq (E+D_z\,K\,C)\,\gS=0_{\scriptscriptstyle \gZ}$, 
so that $T_2(\l)$ is zero.
From (\ref{eqa0}) and $\gS\subseteq \gV$ we find
$(\l\,I-A-B\,K\,C)\,\gS\subseteq \gV$.
Using this with (\ref{u1}) and (\ref{u2}) gives
\beann
&&\hspace{-1cm} (E+D_z F)(\l I-A-B F)^{-1}(\l I-A-B K C)(\l I-A-G C)^{-1}(H+G\,G_y)\\
&& \subseteq 
(E+D_z\,F)(\l\,I-A-B\,F)^{-1}\,(\l\,I-A-B\,K\,C)\,\gS  \subseteq 
(E+D_z\,F)(\l\,I-A-B\,F)^{-1}\,\gV=0_{\scriptscriptstyle \gZ}.
\eeann
Thus, $T_3(\l)$ is zero.
Finally, from (\ref{eqd0}) we find $T_4=G_z+D_z\,K\,G_y={0}$. It follows that $G_{\z,\w}(\l)=0$.

{\bf (Only if)}. Let $A_c$, $B_c$, $C_c$ and $D_c$ exist such that $I-D_y\,D_c$ is non-singular and  $G_{\z,\w}(\l)=0$. This implies that $\widehat{G}=0$, and there exists an $\widehat{A}$-invariant subspace $\widehat{\gI}$ such that $\ima \widehat{H}\subseteq \widehat{\gI}\subseteq \ker \widehat{C}$, see \cite[Thm.~4.6]{Trentelman-SH-01}.
 We start proving that  $\gV=\mathfrak{p}\,(\widehat{\gI})$ is $(A,B,E,D_z)$-output nulling, where $\mathfrak{p}$ denotes the projection on $\gX$ (see Appendix~A). Let $\x\in \gV$. There exists $\p\in \gP$ such that $\bsmat \x\\[1mm] \p\esmat\in \widehat{\gI}$. Since $\widehat{\gI}$ is $\widehat{A}$-invariant, we have $\widehat{A}\,\bsmat \x\\[1mm] \p\esmat\in \widehat{\gI}$, i.e.
\beann
\bmat{c}
A\,\x+B\,D_c\,W\,C\,\x+ B\,C_c\,\p+B\,D_c\,W\,D_y\,C_c\,\p \\
B_c\,W\,C\,\x+A_c\,\p+B_c\,W\,D_y\,C_c\,\p \emat\in \widehat{\gI},
\eeann
which implies
$A\,\x+B\,D_c\,W\,C\,\x+ B\,C_c\,\p+B\,D_c\,W\,D_y\,C_c\,\p \in \mathfrak{p}\,(\widehat{\gI})$.
On the other hand, since $\widehat{\gI} \subseteq \ker \widehat{C}$, we have also
$\widehat{C}\,\bsmat \x\\[1mm] \p\esmat=E\,\x+
D_z\,D_c\,W\,C\,\x+ D_z\,C_c\,\p+D_z\,D_c\,W\,D_y\,C_c\,\p=0_{\scriptscriptstyle \gZ}$.
We can write these two equations together as
\beann
\bmat{c}  A  \\[-1mm]  E  \emat\x+\bmat{c}   B   \\[-1mm]   D_z  \emat\bigl(D_c\,W\,C\,\x+C_c\,\p+D_c\,W\,D_y\,C_c\,\p\bigr) \in \mathfrak{p}\,(\widehat{\gI})\oplus 0_{\scriptscriptstyle \gZ},
\eeann
so that $\bsmat A \\[1mm] E\esmat\,\x\in  \mathfrak{p}\,(\widehat{\gI})\oplus 0_{\scriptscriptstyle \gZ}+\ima \bsmat B \\[1mm] D_z\esmat$. Thus $\gV=\mathfrak{p}\,(\widehat{\gI})$ is $(A,B,E,D_z)$-output nulling
as required. Now we prove that $\gS= \mathfrak{i}\,(\widehat{\gI})$ is $(A,H,C,G_y)$-input containing, where $\mathfrak{i}$ denotes the intersection (see Appendix~A). 
Let $\bsmat \x\\[1mm] \w\esmat\in \gS\oplus \gW \cap \ker [\begin{array}{cc} C & G_y \end{array}]$. We need to prove that
$[\begin{array}{cc} A & H \end{array}]\bsmat \x \\[1mm] \w\esmat\in \gS$.
Since $\x\in \gS= \mathfrak{i}\,(\widehat{\gI})$, we obtain $\bsmat \x\\[1mm] {0}\esmat\in \widehat{\gI}$, and since $\widehat{\gI}$ is $\widehat{A}$-invariant, we find
$\widehat{A}\bsmat \x\\[1mm] 0\esmat=
\bsmat
A\,\x+B\,D_c\,W\,C\,\x \\[1mm]
B_c\,W\,C\,\x \esmat\in \widehat{\gI}$. 
Since $\widehat{\gI}\subseteq \ima \widehat{H}$, we can write $\widehat{H}\,\w\in \widehat{\gI}$, i.e., 
$\bsmat H+B\,D_c\,W\,G_y \\[1mm]
B_c\,W\,G_y \esmat \,\w\in \widehat{\gI}$.
From the last two relations we find
\[
\bmat{cc}
A\,\x+B\,D_c\,W\,C\,\x+H\,\w+B\,D_c\,W\,G_y\,\w \\[0mm]
B_c\,W\,C\,\x+B_c\,W\,G_y\,\w \emat\in \widehat{\gI}.
\]
Since $\bsmat \x \\[1mm] \w \esmat\in \ker [\begin{array}{cc} C & G_y \end{array}]$, the latter can be simplified to
$\bsmat
A\,\x+H\,\w \\[1mm]
{0} \esmat\in \widehat{\gI}$,
i.e., 
$[\begin{array}{cc} A & H \end{array}]\bsmat \x\\[1mm] \w \esmat\in \mathfrak{i}\,(\widehat{\gI})=\gS$, as required.

Now our aim is to show that {\bf (i-iii)} are satisfied.
Since $\widehat{\gI}\supseteq \ima \widehat{H}$, it follows that $\mathfrak{p}\,(\widehat{\gI})\supseteq \mathfrak{p}\,(\ima \widehat{H})$, see Lemma \ref{PI1}, which can be rewritten as
$\gV\supseteq \ima\bigl(H+B\,D_c\,W\,G_y\bigr)$. This inclusion together with $\widehat{G}=0$
leads to $\gV \supseteq \ima(H+B\,\Phi)$ and $G_z+D_z\,{\Phi}=0$, where  $\Phi=D_c\,W\,G_y$. Denoting by $V$ a basis matrix of $\gV$, in view of the these equations there exists a matrix $X$ such that
$H+B\,\Phi=V\,X$ and $G_z+D_z\,\Phi=0$, 
i.e.,
$\bsmat
H\\[1mm] G_z \esmat=\bsmat V \\[1mm] 0\esmat\,X+\bsmat
B\\[1mm] D_z \esmat\,\Phi$,
so that {\bf (i)} is satisfied.
Since $\widehat{C}\,\widehat{\gI}=0_{\scriptscriptstyle \gZ}$, then 
$\bigl(E+D_z\,D_c\,
W\,
C\bigr) \,\mathfrak{i}\,(\widehat{\gI})=0_{\scriptscriptstyle \gZ}$,
see Lemma \ref{PI2}.
Since $\widehat{G}=0$ then
$\bigl(E+{\Psi}\,C\bigr) \,\gS=0_{\scriptscriptstyle \gZ}$ and $G_z+{\Psi}\,G_y=0$.
Let $Q$ be a full row-rank matrix such that $\ker Q=\gS$; we obtain $\ker Q \subseteq \ker \bigl(E+\Psi\,C\bigr)$, so that a matrix $K$ of suitable size exists such that $\Theta\,Q=E+\Psi\,C$. Thus
$E+\Psi\,C=\Theta\,Q$ and
$G_z+\Psi\,G_y=0$, 
i.e., 
$[\begin{array}{cc}  E & G_z \end{array}]=\Theta\, [\begin{array}{cc}  Q & 0 \end{array}]-\Psi \, [\begin{array}{cc}  C & G_y \end{array}]$, which another way of writing
$\ker  [\begin{array}{cc} E & G_z\end{array}]\supseteq \ker [\begin{array}{cc} Q & 0\end{array}]\cap \ker[\begin{array}{cc} C & G_y \end{array}]$.
Since $\ker[\begin{array}{cc} Q & 0 \end{array}]=\gS \oplus \gW$, we obtain
$\ker  [\begin{array}{cc} E & G_z \end{array}]\supseteq (\gS \oplus \gW) \cap \ker[\begin{array}{cc} C & G_y \end{array}]$.
We have proved {\bf (i-ii)}. The proof of {\bf (iii)} follows from $\mathfrak{i}\,(\widehat{\gI})\subseteq \mathfrak{p}\,(\widehat{\gI})$. 

From Lemma~\ref{lem3cond} there exists $K\in \real^{m \times p}$ such that (\ref{clK}) 
 holds. We show that one of such $K$ is also such that $I+K\,D_y$ is non-singular. Let $K=D_c\,W$. From the matrix inversion lemma, $I+K\,D_y$ is non-singular. It remains to prove that $K$ satisfies (\ref{clK}). Rewriting (\ref{clK}) using 
 $K=D_c\,W$ gives
\beann
 \bmat{cc} A+B\,D_c\,W\,C & H+B\,D_c\,W\,G_y \\[-1.5mm]
 E+D_z\,D_c\,W\,C & G_z+ D_z\,D_c\,W\,G_y \emat(\gS \oplus \gW)\subseteq \gV \oplus 0_{\scriptscriptstyle \gZ}.
 \eeann
 Let $\bsmat \v \\[1mm] \w \esmat \in \gS\oplus \gW$. 
 We want to prove that 
 \bea
 \label{toprove}
  \bmat{cc} A\,\v+B\,D_c\,W\,C\,\v+ H\,\w+B\,D_c\,W\,G_y\,\w \\
 E\,\v+D_z\,D_c\,W\,C\,\v+ G_z\,\w+ D_z\,D_c\,W\,G_y\,\w \emat\in \gV \oplus 0_{\scriptscriptstyle \gZ}.
 \eea
 Since $\v\in \gS=\mathfrak{i}(\widehat{\gI})$, we have $\bsmat \v \\[1mm] 0 \esmat\in \widehat{\gI}$. Since $\widehat{\gI}$ is $\widehat{A}$-invariant, we find
$\widehat{A}\,\bsmat \v \\[1mm] 0 \esmat=\bsmat
 A\,\v+B\,D_c\,W\,C\,\v \\[1mm] E\,\v+D_z\,D_c\,W\,C\,\v \esmat\in \widehat{\gI}$.
 It follows that  $A\,\v+B\,D_c\,W\,C\,\v\in \mathfrak{p}(\widehat{\gI})=\gV$. Moreover, since $\ima \widehat{H}\subseteq \widehat{\gI}$, we have 
\[
\bmat{cc} H\,\w+B\,D_c\,W\,G_y\,\w \\[0mm] B_c\,W G_y \,\w \emat \in \widehat{\gI}.
\]
In particular, $H\,\w+B\,D_c\,W\,G_y\,\w\in \mathfrak{p}(\widehat{\gI})=\gV$. We have proved that, in (\ref{toprove}), there holds $A\,\v+B\,D_c\,W\,C\,\v+ H\,\w+B\,D_c\,W\,G_y\,\w\in  \mathfrak{p}(\widehat{\gI})=\gV$. Since the system is disturbance decoupled, the feedthrough $G_z+ D_z\,D_c\,W\,G_y$ is zero. Hence, it remains to show that $E\,\v+D_z\,D_c\,W\,C\,\v=0$. This follows from the fact that $\widehat{C}\,\widehat{\gI}=0$, so that $\widehat{C}\,\bsmat \v \\[1mm] 0 \esmat=0$, which gives $E\,\v+D_z\,D_c\,W\,C\,\v=0$.
\endproof

\begin{remark}
{ The statement of Theorem~\ref{solv} involves conditions that are not independent. Indeed,
Lemma~\ref{lem3cond} showed the relationship between {\bf (i-iii)} and condition (\ref{clK}) in {\bf (iv)}. Thus, if the necessity and the sufficiency statements are kept separate, some of the conditions in the statement of Theorem~\ref{solv} are absorbed into the others. However, we prefer this way of presenting this result, because it displays the symmetry between the two implications of the statement.
}
\end{remark}

\begin{remark}
{
The well-posedness condition on the invertibility of the matrix $I+K\,D_y$ is essential in the nonstrictly proper case. Indeed, 
 there are cases where the entire set of all possible $K$ matrices satisfying (\ref{clK}) renders $I+K\,D_y$ singular. Consider for example
\beann
A &=& \bsmat 
0 && 0 && 0 \\[1mm]
0 && 0 && 0 \\[1mm]
-1 && 0 && 0 \esmat, \;\; B=\bsmat
 0 && 0 \\[1mm]
-1  && 0 \\[1mm]
0 && -1  \esmat, 
\;\; H=\bsmat
1 && 0 \\[1mm]
0  && 1 \\[1mm]
1 && 0  \esmat, 
\\
C&=&\bsmat
-1 && 0 && 0 \\[1mm]
0 && 1 && 1  \esmat, \;\; D_y=\bsmat
1 && 0 \\[1mm]
0 && -1 \esmat, \;\; G_y=\bsmat
0 && 0 \\[1mm]
-1 && -1 \esmat, \\ E&=&\bsmat
0 && 0 && 1 \esmat, \;\; D_z=\bsmat
-1 && 0 \esmat, \;\; G_z=\bsmat
0 && 0 \esmat, \;\; \gS=\spanR\left\{\bsmat 1\\[1mm] -1 \\[1mm] 1 \esmat\right\},
\eeann
and $\gV=\real^3$.
Subspace $\gV$ is $(A,B,E,D_z)$-output nulling and $\gS$ is $(A,H,C,G_y)$-input-input containing, and they satisfy {\bf (i-iii)} of Theorem~\ref{solv}. Thus, a matrix $K$ exists that satisfies (\ref{clK}). 
%
One can easily see that the set of all matrices $K$ for which (\ref{clK}) is fulfilled is given by
$K=
  \bsmat -1 && 0 \\[1mm] \alpha && \beta\esmat$,
 where $\alpha,\beta$ are free parameters. Clearly,
 $I+K\,D_y=\bsmat 0 && 0 \\[1mm] \alpha && -\beta \esmat$, which is singular for every choice of $\alpha,\beta$.
 
}
\end{remark}

\begin{remark}
{The {\em if} part of the proof of Theorem~\ref{solv} offers a compensator structure which involves a feedback matrix $K$ such that (\ref{clK}) is satisfied, an $(A,B,E,D_z)$-ouput-nulling friend $F$ of $\gV$ and an $(A,H,C,G_y)$-input containing friend $G$ of $\gS$. This, however, does not constitute a parameterization of all the decoupling filters. Consider for example a system described by the matrices
\beann
A&=&\bsmat 1 && 0 \\[1mm] 0 && 1 \esmat, \;\; B=\bsmat -1 \\[1mm] 0 \esmat,\;\; H=\bsmat 1\\[1mm] 0 \esmat,\;\;
C=\bsmat 1 && 0 \esmat, \\
\;\; D_y&=&G_y=1, \;\;E=\bsmat 0 && -1 \esmat, \;\; D_z=G_z=0.
\eeann
One can verify that the compensator described by
$A_c=\bsmat 0 && 0 \\[1mm] 0 && 0 \esmat, \;\; B_c=\bsmat 0 \\[1mm] 10 \esmat,\;\; C_c=\bsmat 0 && 3 \esmat, \;\; D_c=6$
solves the disturbance decoupling problem. Inverting the last three equations of (\ref{eqcomp}) we obtain
\beann
K &=& D_c\,(I-D_y\,D_c)^{-1}=-6/5 \\
F &=& (I-D_c\,D_y)^{-1}\,C_c+K\,C=[\begin{array}{cc} -6/5 & -3/5 \end{array}] \\
G &=& (I-D_c\,D_y)^{-1}(B\,D_c-B_c)=\bsmat 6/5 \\[1mm] 2 \esmat.
\eeann
However, when using these values in the first of (\ref{eqcomp}) we obtain 
$A+G\,C+(B+G\,D_y)\,(I+K\,D_y)^{-1}(F-K\,C)=\frac{1}{5} \bsmat 11 && 3 \\[1mm] 10 && 35\esmat$,
which does not coincide with $A_c$.\footnote{Note also that matrix $G$  is not an input containing friend of $\gS$.}
Hence, the decoupling filter proposed here does not fall in the category of those obtainable as in the proof of Theorem~\ref{solv}. Nevertheless, it is still true that a compensator in the desired form can always be found. Indeed, any ($1 \times 1$) matrix $K$ satisfies (\ref{clK}). For example, choosing $K=1/2$ 
and the friends $F=[\begin{array}{cc} 1 & 0 \end{array}]$ and $G=0$, we obtain 
 \beann
 A_{c,1} \ns&\ns = \ns&\ns A+G C+(B+G D_y)(I+K D_y)^{-1}(F-K C)=\bsmat 2/3 && 0 \\[1mm] 0 && 1 \esmat, \\
 B_{c,1} \ns&\ns = \ns&\ns (B+G D_y)(I+K D_y)^{-1} K-G=-\bsmat 1/3\\[1mm] 0 \esmat,\\
 C_{c,1} \ns&\ns = \ns&\ns (I+K D_y)^{-1} (F-K C)=[\begin{array}{cc} \frac{1}{3} & 0 \end{array}],\\
D_{c,1} \ns&\ns = \ns&\ns (I+K D_y)^{-1} \,K={1}/{3}.
\eeann
In other words, if there exists a compensator that solves the decoupling problem, it may not be obtainable in the way described in the proof of Theorem~\ref{solv}. However, we know that we can always find $\gS$ and $\gV$ as the intersection and projection of an invariant for the extended system contained in $\ker \widehat{C}$ and containing $\ima \widehat{H}$ and matrix $K$, and determining the friends of $\gV$ and $\gS$ we can construct an alternative compensator that may not be the one we had originally.
It is now possible to better appreciate the role of condition {\bf (iv)} in Theorem~\ref{solv}, which guarantees that, even  
 if the parameterization of the decoupling filters is not exhaustive, every controller is associated to at least one feasible matrix $K$. 
}
\end{remark}

The solvability conditions of Theorem~\ref{solv} can be also stated in terms of $\gV^\star_{\scriptscriptstyle (A,B,E,D_z)}$ and $\gS^\star_{\scriptscriptstyle (A,H,C,G_y)}$.

\begin{corollary}
\label{corsolv}
Problem \ref{pro1} is solvable if and only if there exist a matrix $K\in \real^{m \times p}$ such that
\begin{description}
\item{\bf (i)} $\ima \bmat{c} H \\[-1.5mm] G_z \emat \subseteq (\gV^\star_{\scriptscriptstyle (A,B,E,D_z)} \oplus 0_{\scriptscriptstyle \gZ})+\ima \bmat{c} B \\[-1.5mm] D_z \emat$;
\item{\bf (ii)} $\ker \,[\begin{array}{cc} E & G_z \end{array}] \supseteq (\gS^\star_{\scriptscriptstyle (A,H,C,G_y)} \oplus \gW)  \cap \ker \,[\begin{array}{cc}  C & G_y \end{array}] $;
\item{\bf (iii)} $\gS^\star_{\scriptscriptstyle (A,H,C,G_y)}\subseteq\gV^\star_{\scriptscriptstyle (A,B,E,D_z)}$;
\item{\bf (iv)} $I+K\,D_y$ is non-singular, and $K$ satisfies  \\[-.7cm]
{\small
\bea
 \label{clK2}
 \bmat{cc} A+B K C & H+B K G_y  \\
 E+D_z K G_y & G_z+ D_z K G_y  \emat(\gS^\star_{\scriptscriptstyle (A,H,C,G_y)} \oplus \gW)\subseteq \gV^\star_{\scriptscriptstyle (A,B,E,D_z)} \oplus 0_{\scriptscriptstyle \gZ}.
\eea
}
\end{description}
\end{corollary}
\proof The sufficiency is obvious from Theorem~\ref{solv}. 
Let us prove the necessity. Let the problem be solvable. In view of Theorem~\ref{solv}, there exist two subspaces $\gV$ and $\gS$ and a matrix $K$ satisfying all the conditions in its statement. We find
\[
\ima \bmat{c} H \\[0mm] G_z \emat \subseteq (\gV \oplus 0_{\scriptscriptstyle \gZ})+\ima \bmat{c} B \\[0mm] D_z \emat  \subseteq (\gV^\star_{\scriptscriptstyle (A,B,E,D_z)} \oplus 0_{\scriptscriptstyle \gZ})+\ima \bmat{c} B \\[0mm] D_z \emat,
\]
 \[
 \ker \,[\begin{array}{cc} E & G_z \end{array}] \supseteq (\gS\oplus \gW)  \cap \ker \,[\begin{array}{cc}  C & G_y \end{array}] \supseteq (\gS^\star_{\scriptscriptstyle (A,H,C,G_y)} \oplus \gW)  \cap \ker \,[\begin{array}{cc}  C & G_y \end{array}],
 \]
  \[
  \gS^\star_{\scriptscriptstyle (A,H,C,G_y)}\subseteq \gS\subseteq \gV \subseteq\gV^\star_{\scriptscriptstyle (A,B,E,D_z)},
  \]
   and
\beann
\bmat{cc} A+B\,K\,C & H+B\,K\,G_y \\[0mm]
 E+D_z\,K\,C & G_z+ D_z\,K\,G_y \emat(\gS^\star_{\scriptscriptstyle (A,H,C,G_y)} \oplus \gW) \ns&\ns \subseteq  \ns&\ns 
\bmat{cc} A+B\,K\,C & H+B\,K\,G_y \\[0mm]
 E+D_z\,K\,C & G_z+ D_z\,K\,G_y \emat(\gS \oplus \gW) \subseteq \gV \oplus 0_{\scriptscriptstyle \gZ}\\
  \ns&\ns \subseteq  \ns&\ns \gV^\star_{\scriptscriptstyle (A,B,E,D_z)} \oplus 0_{\scriptscriptstyle \gZ}.
 \eeann
\endproof

\section{Solution of Problem~\ref{pro2}}
We now consider Problem~\ref{pro2}. Two necessary solvability conditions are the asymptotic stabilizability of the pair $(A,B)$ and the asymptotic detectability of the pair $(C,A)$ \cite[Thm.~3.40]{Trentelman-SH-01}.
These are, therefore, standing assumptions for this section.
The following result provides a solution to Problem~\ref{pro2} in terms of the largest $(A,B,E,D_z)$-stabilizability subspace  and of the smallest $(A,H,C,G_y)$-detectability subspace, see \cite[Thm.~4.1]{Stoorvogel-W-91}.

\begin{theorem}
\label{solvVSg}
Problem~\ref{pro2} is solvable if and only if there exist a matrix $K\in \real^{m \times p}$ such that
\begin{description}
\item{\bf (i)} $\ima \bmat{c} H \\[-1.5mm] G_z \emat \subseteq (\gV^\star_{\scriptscriptstyle (A,B,E,D_z),g} \oplus 0_{\scriptscriptstyle \gZ})+\ima \bmat{c} B \\[-1.5mm] D_z \emat$;
\item{\bf (ii)} $\ker \,[\begin{array}{cc} E & G_z \end{array}] \supseteq (\gS^\star_{\scriptscriptstyle (A,H,C,G_y),g} \oplus \gW)  \cap \ker \,[\begin{array}{cc}  C & G_y \end{array}] $;
\item{\bf (iii)} $\gS^\star_{\scriptscriptstyle (A,H,C,G_y),g}\subseteq\gV^\star_{\scriptscriptstyle (A,B,E,D_z),g}$;
\item{\bf (iv)} $I+K\,D_y$ is non-singular, and $K$ satisfies 
{\small 
\beann
 \bmat{cc}  A+B K C & H+B K G_y  \\
 E+D_z K G_y & G_z+ D_z K G_y  \emat(\gS^\star_{\scriptscriptstyle (A,H,C,G_y),g} \oplus \gW)\subseteq \gV^\star_{\scriptscriptstyle (A,B,E,D_z),g} \oplus 0_{\scriptscriptstyle \gZ}.
\eeann
}
\end{description}
\end{theorem}

An immediate consequence is the following result.

\begin{corollary}
\label{corDDPDOF}
Problem~\ref{pro2} is solvable if and only if there exist an $(A,B,E,D_z)$-stabilizability output nulling subspace $\gV$ and an $(A,H,C,G_y)$-detectability input containing subspace $\gS$ and a matrix $K\in \real^{m \times p}$ such that
\begin{description}
\item{\bf (i)} $\ima \bmat{c} H \\[-1.5mm] G_z \emat\subseteq \gV \oplus {0}_{\gZ} +\ima \bmat{c} B \\[-1.5mm] D_z\emat$;
\item{\bf (ii)} $[\begin{array}{cc} E & G_z \end{array}]\left(\gS \oplus \gW  \cap \ker [\begin{array}{cc} C & G_y \end{array}]\right) =0_{\scriptscriptstyle \gZ}$;
\item{\bf (iii)} $\gS\subseteq \gV$;
\item{\bf (iv)} $I+K\,D_y$ is non-singular, and $K$ satisfies \\[-.7cm]
{\small 
\bea
\label{clK4}
 \bmat{cc}   A+B K C & H+B K G_y  \\
 E+D_z K G_y & G_z+ D_z K G_y  \emat(\gS \oplus \gW)\subseteq \gV \oplus 0_{\scriptscriptstyle \gZ}.
\eea}
\end{description}
\end{corollary}
\proof
{\bf (Only if)}. It follows directly from Theorem~\ref{solvVSg}, by taking $\gV=\gV^\star_{\scriptscriptstyle (A,B,E,D_z),g}$ and $\gS=\gS^\star_{\scriptscriptstyle (A,H,C,G_y),g}$. \\
{\bf (If)}. Since $\gV$ is internally stabilizable, in view of the stabilizability of the pair $(A,B)$, $\gV$ is also externally stabilizable; thus, there exists an output nulling friend $F$ of $\gV$ such that $A+B\,F$ is asymptotically stable. Likewise, since $\gS$ is externally detectable, the detectability of the pair $(C,A)$ ensures that $\gS$ is also internally detectable; it follows that there exists an input containing friend $G$ of $\gS$ such that $A+G\,C$ is asymptotically stable. We can therefore follow the same steps of the proof of 
Theorem~\ref{solv}, and we obtain
that a matrix $K$ exists such that (\ref{clK4}) holds. 
Defining the compensator matrices in the same way as in the proof of Theorem~\ref{solv}, we obtain that the eigenvalues of the closed-loop system are $\sigma(A+B\,F)\uplus \sigma(A+G\,C)$, and $G_{z,w}(\l)$ is  zero.
\endproof

We now generalize the solvability stated in terms of self bounded and self hidden subspaces, namely $\gV_m+\gS_M$ in place of $\gS^\star_{\scriptscriptstyle (A,H,C,G_y)}$ and $\gV_m$ in place of $\gV^\star_{\scriptscriptstyle (A,B,E,D_z)}$. 
The first and more important step, which arises in the nonstrictly proper case, is to prove that the well-posedness condition does not change if we choose these self bounded and self hidden subspaces instead of $\gS^\star_{\scriptscriptstyle (A,H,C,G_y)}$ and  $\gV^\star_{\scriptscriptstyle (A,B,E,D_z)}$. 

\begin{theorem}
\label{thm:K}
Let Problem~\ref{pro1} be solvable.
The set of matrices $K$ that satisfy (\ref{clK2})
coincides with the set of matrices $K$ that satisfy 
{
 \bea
 \label{inc2}
 \bmat{cc}   A+B K C & H+B K G_y  \\
 E+D_z K G_y & G_z+ D_z K G_y   \emat  (\gS_M \oplus \gW)\subseteq (\gS_M+\gV_m) \oplus 0_{\scriptscriptstyle \gZ}.
 \eea
 }
\end{theorem}
%
\proof 
Since $\gS^\star_{\scriptscriptstyle (A,H,C,G_y)}\subseteq \gS_M$ and $\gV^\star_{\scriptscriptstyle (A,B,E,D_z)}\supseteq \gS_M+\gV_m$,
if $K$ satisfies (\ref{inc2}), it also satisfies (\ref{clK2}) since
{
\beann
&& \bmat{cc}  A+B K C & H+B K G_y  \\
 E+D_z K G_y & G_z+ D_z K G_y    \emat(\gS^\star_{\scriptscriptstyle (A,H,C,G_y)}\oplus \gW)  \subseteq 
 \bmat{cc}  A+B K C & H+B K G_y  \\
 E+D_z K G_y & G_z+ D_z K G_y    \emat (\gS_M \oplus \gW)\\
&& \qquad
 \subseteq (\gS_M+\gV_m) \oplus {0}_{\gZ} 
  \subseteq \gV^\star_{\scriptscriptstyle (A,B,E,D_z)} \oplus {0}_{\scriptscriptstyle \gZ}.
  \eeann
  }
  We now prove that if $K$ satisfies (\ref{clK2}), it also satisfies (\ref{inc2}).
  Let $K$ be such that (\ref{clK2}) holds. Proving that $K$ also satisfies (\ref{inc2}) amounts to proving 
the four inclusions
 \bea
&& (A+B\,K\,C)\,\gS_M \subseteq \gS_M +\gV_m, \label{k1} \\
&&  \ima (H+B\,K\,G_y) \subseteq \gS_M +\gV_m, \label{k2} \\
&&  (E+D_z\,K\,G_y)\,\gS_M = 0,  \label{k3}\\
&&  G_z+D_z\,K\,G_z=0. \label{k4}
  \eea
Note that (\ref{k4}) trivially holds because $K$ solves Problem~\ref{pro1} (see proof of Theorem~\ref{solv}).
Consider (\ref{k3}). We show that 
$(E+D_z\,K\,G_y)\,\gS^\star_{\scriptscriptstyle (A,H,C,G_y)}=0_{\scriptscriptstyle \gZ}$ implies $(E+D_z\,K\,G_y)\,\gS_M=0_{\scriptscriptstyle \gZ}$.
Recall that 
   \[
   \gS_M=\gS^\star_{\scriptscriptstyle \left(A,H,\bsmat {\scriptscriptstyle C}\\[0.3mm]{\scriptscriptstyle E} \esmat,\bsmat {\scriptscriptstyle G_y} \\ {\scriptscriptstyle G_z}\esmat\right)}+
   \gV^\star_{\scriptscriptstyle \left(A,H,\bsmat {\scriptscriptstyle C}\\[0.3mm]{\scriptscriptstyle E} \esmat,\bsmat {\scriptscriptstyle G_y} \\ {\scriptscriptstyle G_z}\esmat\right)}
   =\gS^\star_{\scriptscriptstyle (A,H,C,G_y)}+
   \gV^\star_{\scriptscriptstyle \left(A,H,\bsmat {\scriptscriptstyle C}\\[0.3mm]{\scriptscriptstyle E} \esmat,\bsmat {\scriptscriptstyle G_y} \\ {\scriptscriptstyle G_z}\esmat\right)},
   \]
    where the last equality follows from the fact that Problem~\ref{pro1} is solved.
  From $(E+D_z\,K\,G_y)\,\gS^\star_{\scriptscriptstyle (A,H,C,G_y)}=0$ and Lemma~\ref{lem41} 
   \beann
    (E+D_z\,K\,G_y)\,\gS_M =(E+D_z\,K\,G_y)\,\gV^\star_{\scriptscriptstyle \left(A,H,\bsmat {\scriptscriptstyle C}\\[0.3mm]{\scriptscriptstyle E} \esmat,\bsmat {\scriptscriptstyle G_y} \\ {\scriptscriptstyle G_z}\esmat\right)}
    \subseteq (E+D_z\,K\,G_y) \left(\bmat{c} C\\[-2mm] E \emat^{-1} \ima \bmat{c} G_y \\[-2mm] G_z \emat\right).
   \eeann
   We prove that 
   \[
   (E+D_z\,K\,G_y) \left(\bmat{c} C\\[0mm] E \emat^{-1} \ima \bmat{c} G_y \\[0mm] G_z \emat\right)=0,
   \]
    i.e., $\bsmat C\\[1mm]E \esmat^{-1} \ima \bsmat G_y \\[1mm] G_z \esmat\subseteq \ker (E+D_z\,K\,G_y)$. 
Let $\x$ be a vector of the left hand-side, so that 
$\bsmat C\\[1mm] E \esmat\,\x\in \ima \bsmat G_y \\[1mm] G_z \esmat$,
   so that there exists $\w$ such that $C\,\x=G_y\,\w$ and $E\,\x=G_z\,\w$. Thus,
   \[
   (E+D_z\,K\,G_y)\,\x=G_z\,\w+D_z\,K\,G_y\,\w=(G_z+D_z\,K\,G_z)\,\w=0,
   \]
as required. Consider (\ref{k2}). We need to prove that 
$\ima (H+B\,K\,G_y) \subseteq \gV^\star_{\scriptscriptstyle (A,B,E,D_z)}$ implies $\ima (H+B\,K\,G_y) \subseteq \gS_M +\gV_m$. 
 Using the last inclusion $ \gS_M +\gV_m\supseteq\gV^\star_{\scriptscriptstyle (A,B,E,D_z)}  \cap [\begin{array}{cc} B & H \end{array}]\,\ker[\begin{array}{cc} D_z & G_z \end{array}]$ in the proof of Corollary~\ref{cor5}, we only need to prove that 
$\ima (H+B\,K\,G_y) \subseteq[\begin{array}{cc} B & H \end{array}]\,\ker[\begin{array}{cc} D_z & G_z \end{array}]$.
 Let $\x\in  \ima (H+B\,K\,G_y)$. There exists $\w$ such that 
 $\x= (H+B\,K\,G_y)\,\w$. Let $g=K\,G_y\,\w$, so that $\x=H\,\w+B\,g$ and from {\em (iv)} we also have $G_z+D_z\,K\,G_z=0$. Multiplying this by $\w$ gives 
 $G_z\,\w+D_z\,g=0$. Therefore, 
$\x=[\begin{array}{cc} B & H \end{array}]\,\bsmat \g \\[1mm] \w \esmat$, where $G_z\,\w+D_z\,g=0$. Thus $\x\in [\begin{array}{cc} B & H \end{array}]\,\ker[\begin{array}{cc} D_z & G_z \end{array}]$ as required.
We now prove (\ref{k1}). We have to prove that 
$(A+B\,K\,C)\,\gS^\star_{\scriptscriptstyle (A,H,C,G_y)}\subseteq \gV^\star_{\scriptscriptstyle (A,B,E,D_z)}$ implies $(A+B\,K\,C)\,\gS_M \subseteq \gS_M +\gV_m$. 
 Recall again that, since Problem~\ref{pro1} is solved, 
\[
\gS_M = \gV^\star_{\scriptscriptstyle \left(A,H,\bsmat {\scriptscriptstyle C}\\[0.3mm]{\scriptscriptstyle E} \esmat,\bsmat {\scriptscriptstyle G_y} \\ {\scriptscriptstyle G_z}\esmat\right)}+\gS^\star_{\scriptscriptstyle \left(A,H,\bsmat {\scriptscriptstyle C}\\[0.3mm]{\scriptscriptstyle E} \esmat,\bsmat {\scriptscriptstyle G_y} \\ {\scriptscriptstyle G_z}\esmat\right)}=
 \gV^\star_{\scriptscriptstyle \left(A,H,\bsmat {\scriptscriptstyle C}\\[0.3mm]{\scriptscriptstyle E} \esmat,\bsmat {\scriptscriptstyle G_y} \\ {\scriptscriptstyle G_z}\esmat\right)}+\gS^\star_{\scriptscriptstyle (A,H,C,G_y)}
 \]
 and
 \[
 \gS_M+\gV_m = \gV_m+\gV^\star_{\scriptscriptstyle \left(A,H,\bsmat {\scriptscriptstyle C}\\[0.3mm]{\scriptscriptstyle E} \esmat,\bsmat {\scriptscriptstyle G_y} \\ {\scriptscriptstyle G_z}\esmat\right)}
 \]
 from the proof of Lemma~\ref{lem22}. Thus
 \beann
 \begin{array}{rcl}
&& \hspace{-4mm} \gS_M+\gV_m \\
  \ns&\ns \ns&\ns =\left(\gS^\star_{\scriptscriptstyle \left(A,[\,B\;\;H\,],E,[\,D_z\;\;G_z\,]\right)} \cap
  \gV^\star_{\scriptscriptstyle \left(A,[\,B\;\;H\,],E,[\,D_z\;\;G_z\,]\right)}\right)+\gV^\star_{\scriptscriptstyle \left(A,H,\bsmat {\scriptscriptstyle C}\\[0.3mm]{\scriptscriptstyle E} \esmat,\bsmat {\scriptscriptstyle G_y} \\ {\scriptscriptstyle G_z}\esmat\right)} \\
  \ns&\ns  \ns&\ns =(\gS^\star_{\scriptscriptstyle \left(A,[\,B\;\;H\,],E,[\,D_z\;\;G_z\,]\right)}+\gV^\star_{\scriptscriptstyle \left(A,H,\bsmat {\scriptscriptstyle C}\\[0.3mm]{\scriptscriptstyle E} \esmat,\bsmat {\scriptscriptstyle G_y} \\ {\scriptscriptstyle G_z}\esmat\right)})\cap \gV^\star_{\scriptscriptstyle \left(A,[\,B\;\;H\,],E,[\,D_z\;\;G_z\,]\right)} \\
  \ns&\ns  \ns&\ns =(\gS^\star_{\scriptscriptstyle \left(A,[\,B\;\;H\,],E,[\,D_z\;\;G_z\,]\right)}+\gV^\star_{\scriptscriptstyle \left(A,H,\bsmat {\scriptscriptstyle C}\\[0.3mm]{\scriptscriptstyle E} \esmat,\bsmat {\scriptscriptstyle G_y} \\ {\scriptscriptstyle G_z}\esmat\right)})\cap \gV^\star_{\scriptscriptstyle (A,B,E,D_z)}.\end{array}
  \eeann
 Using these, we need to show that 
 \beann
&& (A+B\,K\,C)\,(\gV^\star_{\scriptscriptstyle \left(A,H,\bsmat {\scriptscriptstyle C}\\[0.3mm]{\scriptscriptstyle E} \esmat,\bsmat {\scriptscriptstyle G_y} \\ {\scriptscriptstyle G_z}\esmat\right)}+\gS^\star_{\scriptscriptstyle \left(A,H,C,G_y\right)}) \\
&&\qquad  \subseteq (\gS^\star_{\scriptscriptstyle \left(A,[\,B\;\;H\,],E,[\,D_z\;\;G_z\,]\right)}+\gV^\star_{\scriptscriptstyle \left(A,H,\bsmat {\scriptscriptstyle C}\\[0.3mm]{\scriptscriptstyle E} \esmat,\bsmat {\scriptscriptstyle G_y} \\ {\scriptscriptstyle G_z}\esmat\right)})\cap \gV^\star_{\scriptscriptstyle (A,B,E,D_z)}.
\eeann
 This reduces to the four inclusions
 \begin{description}
 \item{\hspace{-4mm}(a)} $(A+B K C)\,\gV^\star_{\scriptscriptstyle \left(A,H,\bsmat {\scriptscriptstyle C}\\[0.3mm]{\scriptscriptstyle E} \esmat,\bsmat {\scriptscriptstyle G_y} \\ {\scriptscriptstyle G_z}\esmat\right)}
 \subseteq \gS^\star_{\scriptscriptstyle (A,[\,B\;\;H\,],E,[\,D_z\;\;G_z\,])}+\gV^\star_{\scriptscriptstyle \left(A,H,\bsmat {\scriptscriptstyle C}\\[0.3mm]{\scriptscriptstyle E} \esmat,\bsmat {\scriptscriptstyle G_y} \\ {\scriptscriptstyle G_z}\esmat\right)}$
  \item{\hspace{-4mm}(b)} $(A+B K C)\,\gV^\star_{\scriptscriptstyle \left(A,H,\bsmat {\scriptscriptstyle C}\\[0.3mm]{\scriptscriptstyle E} \esmat,\bsmat {\scriptscriptstyle G_y} \\ {\scriptscriptstyle G_z}\esmat\right)}
 \subseteq\gV^\star_{\scriptscriptstyle (A,B,E,D_z)}$
 \item{\hspace{-4mm}(c)} $(A+B K C)\,\gS^\star_{\scriptscriptstyle (A,H,C,G_y)}\subseteq  \gS^\star_{\scriptscriptstyle (A,[\,B\;\;H\,],E,[\,D_z\;\;G_z\,])}+\gV^\star_{\scriptscriptstyle \left(A,H,\bsmat {\scriptscriptstyle C}\\[0.3mm]{\scriptscriptstyle E} \esmat,\bsmat {\scriptscriptstyle G_y} \\ {\scriptscriptstyle G_z}\esmat\right)}$
  \item{\hspace{-4mm}(d)} $(A+B K C)\,\gS^\star_{\scriptscriptstyle (A,H,C,G_y)}\subseteq\gV^\star_{\scriptscriptstyle (A,B,E,D_z)}$.
  \end{description}
  Clearly (d) is satisfied because Problem~\ref{pro1} is solvable. We prove (b). The subspace $\gV^\star_{\scriptscriptstyle \left(A,H,\bsmat {\scriptscriptstyle C}\\[0.3mm]{\scriptscriptstyle E} \esmat,\bsmat {\scriptscriptstyle G_y} \\ {\scriptscriptstyle G_z}\esmat\right)}$ satisfies
  \[
    \bmat{c}  A  \\[-1.5mm]   C  \\[-1.5mm]   E   \emat \gV^\star_{\scriptscriptstyle \left(A,H,\bsmat {\scriptscriptstyle C}\\[0.3mm]{\scriptscriptstyle E} \esmat,\bsmat {\scriptscriptstyle G_y} \\ {\scriptscriptstyle G_z}\esmat\right)} \subseteq \left(\gV^\star_{\scriptscriptstyle \left(A,H,\bsmat {\scriptscriptstyle C}\\[0.3mm]{\scriptscriptstyle E} \esmat,\bsmat {\scriptscriptstyle G_y} \\ {\scriptscriptstyle G_z}\esmat\right)}\oplus 0_{\scriptscriptstyle \gY} \oplus 0_{\scriptscriptstyle \gZ}\right) +\ima \bmat{c}   H  \\[-1.5mm]   G_y   \\[-1.5mm]   G_z   \emat.
  \]  
  Let $\tilde{V}$ be a basis matrix of $\gV^\star_{\scriptscriptstyle \left(A,H,\bsmat {\scriptscriptstyle C}\\[0.3mm]{\scriptscriptstyle E} \esmat,\bsmat {\scriptscriptstyle G_y} \\ {\scriptscriptstyle G_z}\esmat\right)}$, so the latter inclusion ensures in particular the existence of matrices $\Xi$ and $\Theta$ of suitable sizes such that $C\,\tilde{V}=G_y\,\Theta$ and $A\,\tilde{V}=\tilde{V}\,X+H\,\Theta$. It follows that
$(A+B\,K\,C)\,\tilde{V}=\tilde{V}\,X+(H+B\,K\,G_y)\,\Theta$,
  so that 
  \[
    (A+B\,K\,C)\,\gV^\star_{\scriptscriptstyle \left(A,H,\bsmat {\scriptscriptstyle C}\\[0.3mm]{\scriptscriptstyle E} \esmat,\bsmat {\scriptscriptstyle G_y} \\ {\scriptscriptstyle G_z}\esmat\right)} \subseteq \gV^\star_{\scriptscriptstyle \left(A,H,\bsmat {\scriptscriptstyle C}\\[0.3mm]{\scriptscriptstyle E} \esmat,\bsmat {\scriptscriptstyle G_y} \\ {\scriptscriptstyle G_z}\esmat\right)} +\ima (H+B\,K\,G_y).
    \]
The inclusion $\gV^\star_{\scriptscriptstyle \left(A,H,\bsmat {\scriptscriptstyle C}\\[0.3mm]{\scriptscriptstyle E} \esmat,\bsmat {\scriptscriptstyle G_y} \\ {\scriptscriptstyle G_z}\esmat\right)} \subseteq\gV^\star_{\scriptscriptstyle (A,B,E,D_z)}$ holds from Lemma~\ref{gt}, together with $\ima (H+B\,K\,G_y) \subseteq \gV^\star_{\scriptscriptstyle (A,B,E,D_z)}$, so that
$(A+B\,K\,C)\,\gV^\star_{\scriptscriptstyle \left(A,H,\bsmat {\scriptscriptstyle C}\\[0.3mm]{\scriptscriptstyle E} \esmat,\bsmat {\scriptscriptstyle G_y} \\ {\scriptscriptstyle G_z}\esmat\right)} \subseteq \gV^\star_{\scriptscriptstyle (A,B,E,D_z)}$.
   Now we prove (a). We have already shown that $\ima (H+B\,K\,G_y) \subseteq[\begin{array}{cc} B & H \end{array}]\,\ker[\begin{array}{cc} D_z & G_z \end{array}]$. 
   Since $\gS^\star_{\scriptscriptstyle (A,[\,B\;\;H\,],E,[\,D_z\;\;G_z\,])}\supseteq[\begin{array}{cc} B & H \end{array}]\,\ker[\begin{array}{cc} D_z & G_z \end{array}]$ we have
$\ima (H+B\,K\,G_y) \subseteq \gS^\star_{\scriptscriptstyle (A,[\,B\;\;H\,],E,[\,D_z\;\;G_z\,])}$.
Adding to both member of this inclusion the subspace $\gV^\star_{\scriptscriptstyle \left(A,H,\bsmat {\scriptscriptstyle C}\\[0.3mm]{\scriptscriptstyle E} \esmat,\bsmat {\scriptscriptstyle G_y} \\ {\scriptscriptstyle G_z}\esmat\right)}$ gives
\[
\gV^\star_{\scriptscriptstyle \left(A,H,\bsmat {\scriptscriptstyle C}\\[0.3mm]{\scriptscriptstyle E} \esmat,\bsmat {\scriptscriptstyle G_y} \\ {\scriptscriptstyle G_z}\esmat\right)}+\ima (H+B K G_y) \subseteq \gV^\star_{\scriptscriptstyle \left(A,H,\bsmat {\scriptscriptstyle C}\\[0.3mm]{\scriptscriptstyle E} \esmat,\bsmat {\scriptscriptstyle G_y} \\ {\scriptscriptstyle G_z}\esmat\right)}+\gS^\star_{\scriptscriptstyle (A,[\,B\;\;H\,],E,[\,D_z\;\;G_z\,])}.
\] 
We have also shown that
\[
(A+B\,K\,C)\,\gV^\star_{\scriptscriptstyle \left(A,H,\bsmat {\scriptscriptstyle C}\\[0.3mm]{\scriptscriptstyle E} \esmat,\bsmat {\scriptscriptstyle G_y} \\ {\scriptscriptstyle G_z}\esmat\right)} \subseteq \gV^\star_{\scriptscriptstyle \left(A,H,\bsmat {\scriptscriptstyle C}\\[0.3mm]{\scriptscriptstyle E} \esmat,\bsmat {\scriptscriptstyle G_y} \\ {\scriptscriptstyle G_z}\esmat\right)}+\ima (H+B\,K\,G_y),
\] 
which readily gives
\[
(A+B\,K\,C)\,\gV^\star_{\scriptscriptstyle \left(A,H,\bsmat {\scriptscriptstyle C}\\[0.3mm]{\scriptscriptstyle E} \esmat,\bsmat {\scriptscriptstyle G_y} \\ {\scriptscriptstyle G_z}\esmat\right)}\subseteq \gV^\star_{\scriptscriptstyle \left(A,H,\bsmat {\scriptscriptstyle C}\\[0.3mm]{\scriptscriptstyle E} \esmat,\bsmat {\scriptscriptstyle G_y} \\ {\scriptscriptstyle G_z}\esmat\right)}+\gS^\star_{\scriptscriptstyle (A,[\,B\;\;H\,],E,[\,D_z\;\;G_z\,])}.
\]
Finally, we prove (c). 
We show in particular that
\bea
\label{eqdmod}
(A+B\,K\,C)\,\gS^\star_{\scriptscriptstyle (A,H,C,G_y)}\subseteq\gS^\star_{\scriptscriptstyle (A,[\,B\;\;H\,],E,[\,D_z\;\;G_z\,])}.
\eea
The inclusion (b) can be written as\footnote{Recall that, if $A\in \real^{m \times n}$, $\gU$ is a subspace of $\real^n$ and $\gH$ is a subspace of $\real^m$, then $A\,\gU\subseteq \gH$ is equivalent to $A^\top \gH^\perp\subseteq \gU^\perp$.} 
\bea
(A^\top+C^\top K^\top B^\top) \gS^\star_{\scriptscriptstyle (A^\top,E^\top,B^\top,D_z^\top)}\subseteq 
\gS^\star_{\scriptscriptstyle (A^\top,[\,C^\top\;\;E^\top\,],H^\top,[\,G_y^\top\;\;G_z^\top\,])}.\label{eqbdual}
\eea
and (\ref{eqbdual}) is equivalent to (\ref{eqdmod}).
From (\ref{eqdmod}), it is trivial to see that (b) holds as well.
\endproof

\begin{theorem}
Problem~\ref{pro2} is solvable  if and only if there exist a matrix $K\in \real^{m \times p}$ such that
\begin{description}
\item{\bf (A)} $\ima \bmat{c} H \\[-1.5mm] G_z \emat\subseteq (\gV^\star_{\scriptscriptstyle (A,B,E,D_z)} \oplus 0_{\scriptscriptstyle \gZ}) +\ima \bmat{c} B \\[-1.5mm] D_z\emat$;
\item{\bf (B)} $\ker [\begin{array}{cc} E & G_z \end{array}]\supseteq \left((\gS^\star_{\scriptscriptstyle (A,H,C,G_y)} \oplus \gW)  \cap \ker [\begin{array}{cc} C & G_y \end{array}]\right)$;
\item{\bf (C)} $\gS^\star_{\scriptscriptstyle (A,H,C,G_y)} \subseteq \gV^\star_{\scriptscriptstyle (A,B,E,D_z)}$;
\item{\bf (D)} $\gV_m+\gS_M$ is an internally stabilizable $(A,B,E,D_z)$-output nulling subspace;
\item{\bf (E)} $\gS_M$ is an externally detectable $(A,H,C,G_y)$-input containing subspace; 
\item{\bf (F)} $I+K\,D_y$ is non-singular, and $K$ satisfies\\[-7mm]
{\small 
\bea
\label{clK11}
 \bmat{cc}  A+B K C & H+B K G_y  \\[-1.5mm]
 E+D_z K G_y & G_z+ D_z K G_y  \emat(\gS^\star_{\scriptscriptstyle (A,H,C,G_y)}  \oplus \gW)\subseteq \gV^\star_{\scriptscriptstyle (A,B,E,D_z)} \oplus 0_{\scriptscriptstyle \gZ}.
\eea
}
\end{description}
\end{theorem}
\proof {\bf (If)}
In view of Corollary~\ref{cor5}, the subspace $\gV_m+\gS_M$ is $(A,B,E,D_z)$-output nulling, while $\gS_M$ is obviously $(A,H,C,G_y)$-input containing. 
Since, from {\bf (D)}-{\bf (E)}, $\gV_m+\gS_M$ is an internally stabilizable $(A,B,E,D_z)$-output nulling
subspace 
and $\gS_M$ is an externally detectable $(A,H,C,G_y)$-input containing subspace,
we can chose as $(A,B,E,D_z)$-stabilizability output nulling subspace the subspace $\gV=\gV_m+\gS_M$ and as $(A,H,C,G_y)$-detectability input containing subspace the subspace $\gS=\gS_M$. We show that the condition of Corollary~\ref{corDDPDOF} are satisfied with this choice of $\gS$ and $\gV$.  
Condition {\bf (iii)} is true by construction. Theorem~\ref{thm:K} guarantees that condition {\bf (iv)} is also satisfied. Finally, in view of Lemma~\ref{lem3cond}, the existence of a matrix $K$ satisfying {\bf (iv)} implies that also conditions {\bf (i)} and {\bf (ii)} hold.\\
{\bf (Only if)}. We assume that Problem~\ref{pro2} is solvable. In view of 
Corollary~\ref{corDDPDOF}, there exist an $(A,B,E,D_z)$-stabilizability output nulling subspace $\gV$ and an $(A,H,C,G_y)$-detectability input containing subspace $\gS$ such that conditions {\bf (i-iv)} in 
Corollary~\ref{corDDPDOF} hold.
Since $\gS^\star_{\scriptscriptstyle (A,H,C,G_y)}\subseteq \gS$ (minimality), $\gV \subseteq \gV^\star_{\scriptscriptstyle (A,B,E,D_z)}$ (maximality), and $\gS\subseteq\gV$, we find that {\bf (A)},{\bf (B)}, {\bf (C)} and {\bf (F)} are satisfied. Now we prove {\bf (D)} and {\bf (E)}. To this end, we show that there exists an internally stabilizable $(A,B,E,D_z)$-self bounded subspace $\bar{\gV}$
such that 
$\gV_m\subseteq \bar{\gV}\subseteq \gV_m+\gS_M$
and an externally detectable $(A,H,C,G_y)$-self hidden subspace $\bar{\gS}$ such that 
$\gV_m\cap \gS_M\subseteq \bar{\gS}\subseteq \gS_M$.
Indeed, consider 
\beann
\begin{array}{rcl}
\bar{\gV} \ns&\ns \defi \ns&\ns \bigl(\gV\cap (\gV_m+\gS_M)\bigr)+\gV_m = (\gV+\gV_m)\cap (\gV_m+\gS_M),  \\
\bar{\gS} \ns&\ns \defi \ns&\ns \bigl(\gS + (\gV_m \cap \gS_M)\bigr)\cap \gS_M = (\gS\cap \gS_M)+(\gV_m \cap \gS_M).\end{array}
\eeann
Since $\ima \bsmat H \\[1mm] G_z \esmat \subseteq (\gV^\star_{\scriptscriptstyle (A,B,E,D_z)} \oplus 0_{\scriptscriptstyle \gZ})+\ima \bsmat B\\[1mm] D_z \esmat$, 
from Corollary~\ref{cor5} the subspace $\gV_m+\gS_M$ is $(A,B,E,D_z)$-self bounded. Moreover, since $\ker \,[\begin{array}{cc} E & G_z \end{array}] \supseteq (\gS^\star_{\scriptscriptstyle (A,H,C,G_y)}  \oplus \gW)  \cap \ker \,[\begin{array}{cc}  C & G_y \end{array}]$, then  $\gV_m\cap\gS_M$  is $(A,B,C,G_y)$-self hidden.
From \cite[Thm.~2]{N-TAC-08}, $\gV_m$ is an internally stabilizable $(A,B,E,D_z)$-output nulling subspace, and, dually, $\gS_M$ is an externally detectable $(A,H,C,G_y)$-input containing subspace, so that {\bf (i)} is proved. Since both $\gV$ and $\gV_m$ are $(A,B,E,D_z)$-output nulling subspaces, so is also their sum $\gV+\gV_m$. Moreover, $\gV+\gV_m$ is also $(A,B,E,D_z)$-self bounded; it follows that the intersection $\bar{\gV}=(\gV+\gV_m)\cap (\gV_m+\gS_M)$ is $(A,B,E,D_z)$-self bounded. Since $\gV+\gV_m$ is $(A,B,E,D_z)$-self bounded and contains $\bar{\gV}$, which is also $(A,B,E,D_z)$-self bounded, an output nulling friend $F$ of 
$\gV+\gV_m$ is also an output nulling friend of $\bar{\gV}$. Since we can choose $F$ so that 
$\gV+\gV_m$ is internally stabilized, the same $F$ stabilizes $\bar{\gV}$ internally, i.e., $\bar{\gV}$ is an internally stabilizable $(A,B,E,D_z)$-output nulling subspace.
Dually, since both $\gS$ and $\gS_M$ are externally detectable $(A,H,C,G_y)$-input containing subspaces, their intersection $\gS\cap\gS_M$ is also externally detectable. Moreover, $\gS\cap\gS_M$ is also $(A,H,C,G_y)$-self hidden; thus, their sum $\bar{\gS}=(\gS\cap \gS_M)+(\gV_m \cap \gS_M)$ is $(A,H,C,G_y)$-self hidden. Since $\gS\cap\gS_M$ is $(A,H,C,G_y)$-self hidden and contained in $\bar{\gS}$, which is also $(A,H,C,G_y)$-self hidden, an input containing friend $G$ of 
$\gS\cap\gS_M$ is also an input containing friend of $\bar{\gS}$. Since we can choose $G$ so that 
$\gS\cap\gS_M$ is externally detected, the same $G$ renders $\bar{\gS}$ detected externally, so that $\bar{\gS}$ is externally detectable.

From Theorem~\ref{the14}, $\ima \bsmat H \\[1mm] G_z \esmat \subseteq (\gV^\star_{\scriptscriptstyle (A,B,E,D_z)} \oplus 0_{\scriptscriptstyle \gZ})+\ima \bsmat B\\[1mm] D_z \esmat$ implies $\ima \bsmat H \\[1mm] G_z \esmat \subseteq (\gV_m \oplus 0_{\scriptscriptstyle \gZ})+\ima \bsmat B\\[1mm] D_z \esmat$, and from its dual  $\ker \,[\begin{array}{cc} E & G_z \end{array}] \supseteq (\gS^\star_{\scriptscriptstyle (A,H,C,G_y)} \oplus \gW)  \cap \ker \,[\begin{array}{cc}  C & G_y \end{array}]$ implies $\ker \,[\begin{array}{cc} E & G_z \end{array}] \supseteq (\gS_M \oplus \gW)  \cap \ker \,[\begin{array}{cc}  C & G_y \end{array}]$.
It follows that
$\ima \bsmat H \\[1mm] G_z \esmat\subseteq (\gV_m \oplus 0_{\scriptscriptstyle \gZ}) +\ima \bsmat B \\[1mm] D_z\esmat
\subseteq (\bar{\gV} \oplus 0_{\scriptscriptstyle \gZ}) +\ima \bsmat B \\[1mm] D_z\esmat$
and $\ker \,[\begin{array}{cc} E & G_z \end{array}] \supseteq (\gS_M \oplus \gW)  \cap \ker \,[\begin{array}{cc}  C & G_y \end{array}]\supseteq (\bar{\gS} \oplus \gW)  \cap \ker \,[\begin{array}{cc}  C & G_y \end{array}]$.
 Finally, from $\gS\subseteq \gV$ we also have the following obvious inclusions
 \beann
 && \gS\cap \gS_M \subseteq \gS\subseteq \gV\subseteq \gV+\gV_m, \\
  && \gS\cap \gS_M \subseteq \gS_M\subseteq \gS_M+\gV_m,
  \eeann
  which imply that $\gS\cap \gS_M$ is contained in the intersection $(\gV+\gV_m)\cap (\gS_M+\gV_m)=\bar{\gV}$; likewise
\beann
  && \gS_M\cap \gV_m \subseteq \gV_m\subseteq \gV+\gV_m \\
  && \gS_M\cap \gV_m \subseteq \gV_m\subseteq \gV_m+\gS_M
  \eeann 
  imply that $\gS_M\cap \gV_m$ is contained in the intersection $\bar{\gV}=(\gV+\gV_m)\cap (\gS_M+\gV_m)$. Their sum $\bar{\gS}=(\gS\cap \gS_M)+(\gS_M\cap \gV_m)$ is therefore also contained in $\bar{\gV}$. Thus, $\bar{\gS}\subseteq \bar{\gV}$.

We already observed that $\gS_M$ is externally detectable.
We now prove that $\gV_m+\gS_M$ is internally stabilizable.
We use the change of coordinate given by a matrix $T=[\begin{array}{cccc} T_1 & T_2 & T_3 & T_4 \end{array}]$ such that
$\ima T_1 = \gV_m\cap \gS_M$, $\ima [\begin{array}{cccc} T_1 & T_2  \end{array}]=\gV_m$, 
$\ima [\begin{array}{cccc} T_1 & T_3 \end{array}]= \gS_M$, $\ima [\begin{array}{cccc} T_1 & T_2 & T_3 \end{array}]= \gS_M+\gV_m$.
We now show that it is always possible to choose $T_3$ in such a way that $\ima T_3\subseteq C^{-1}\ima G_y$. To this end, we prove that 
$\ima [\begin{array}{cccc} T_1 & T_3 \end{array}]=C^{-1}\ima G_y+\ima T_1$, which means that it is always possible to choose $T_3$ in such a way that $\ima T_3\subseteq C^{-1}\ima G_y$.
We have by definition
\beann 
\begin{array}{rcl}
\gS_M \ns&\ns = \ns&\ns \gQ^\star_{\scriptscriptstyle \left(A,H,\bsmat {\scriptscriptstyle C}\\[0.3mm]{\scriptscriptstyle E} \esmat,\bsmat {\scriptscriptstyle G_y} \\ {\scriptscriptstyle G_z}\esmat\right)} \\
\ns&\ns = \ns&\ns \gV^\star_{\scriptscriptstyle \left(A,H,\bsmat {\scriptscriptstyle C}\\[0.3mm]{\scriptscriptstyle E} \esmat,\bsmat {\scriptscriptstyle G_y} \\ {\scriptscriptstyle G_z}\esmat\right)}+\gS^\star_{\scriptscriptstyle \left(A,H,\bsmat {\scriptscriptstyle C}\\[0.3mm]{\scriptscriptstyle E} \esmat,\bsmat {\scriptscriptstyle G_y} \\ {\scriptscriptstyle G_z}\esmat\right)} \\
\ns&\ns = \ns&\ns \gV^\star_{\scriptscriptstyle \left(A,H,\bsmat {\scriptscriptstyle C}\\[0.3mm]{\scriptscriptstyle E} \esmat,\bsmat {\scriptscriptstyle G_y} \\ {\scriptscriptstyle G_z}\esmat\right)}+\gS^\star_{\scriptscriptstyle \left(A,H,C,G_y\right)},
\end{array}
\eeann
where the equality $\gS^\star_{\scriptscriptstyle \left(A,H,\bsmat {\scriptscriptstyle C}\\[0.3mm]{\scriptscriptstyle E} \esmat,\bsmat {\scriptscriptstyle G_y} \\ {\scriptscriptstyle G_z}\esmat\right)} =\gS^\star_{\scriptscriptstyle \left(A,H,C,G_y\right)}$ is a consequence of Theorem~\ref{the13d}.
In view of Lemma~\ref{lem41} we have $\gV^\star_{\scriptscriptstyle \left(A,H,\bsmat {\scriptscriptstyle C}\\[0.3mm]{\scriptscriptstyle E} \esmat,\bsmat {\scriptscriptstyle G_y} \\ {\scriptscriptstyle G_z}\esmat\right)}\subseteq\bsmat C\\[1mm] E\esmat^{-1}\ima\bsmat G_y\\[1mm]G_z\esmat\subseteq C^{-1}\ima G_y$, which implies
$\gS_M\subseteq C^{-1}\ima G_y+\gS^\star_{\scriptscriptstyle \left(A,H,C,G_y\right)}$.
We find
\bea
\gV_m\cap \gS_M
\ns&\ns = \ns&\ns \gR^\star_{\scriptscriptstyle (A,[\,B\;\;H\,],E,[\,D_z\;\;G_z\,]} \cap
\bigl(\gV^\star_{\scriptscriptstyle \left(A,H,\bsmat {\scriptscriptstyle C}\\[0.3mm]{\scriptscriptstyle E} \esmat,\bsmat {\scriptscriptstyle G_y} \\ {\scriptscriptstyle G_z}\esmat\right)}+\gS^\star_{\scriptscriptstyle \left(A,H,\bsmat {\scriptscriptstyle C}\\[0.3mm]{\scriptscriptstyle E} \esmat,\bsmat {\scriptscriptstyle G_y} \\ {\scriptscriptstyle G_z}\esmat\right)}\bigr) \nonumber \\
\ns&\ns = \ns&\ns
\left(\gV^\star_{\scriptscriptstyle (A,[\,B\;\;H\,],E,[\,D_z\;\;G_z\,]}\cap \gS^\star_{\scriptscriptstyle (A,[\,B\;\;H\,],E,[\,D_z\;\;G_z\,]}\right) 
\cap
\Big(\gV^\star_{\scriptscriptstyle \left(A,H,\bsmat {\scriptscriptstyle C}\\[0.3mm]{\scriptscriptstyle E} \esmat,\bsmat {\scriptscriptstyle G_y} \\ {\scriptscriptstyle G_z}\esmat\right)}+\gS^\star_{\scriptscriptstyle (A,H,C,G_y)} \Big)\nonumber  \\
\ns&\ns = \ns&\ns
\left(\gV^\star_{\scriptscriptstyle (A,B,E,D_z)}\cap \gS^\star_{\scriptscriptstyle (A,[\,B\;\;H\,],E,[\,D_z\;\;G_z\,]}\right) 
\cap
\Big(\gV^\star_{\scriptscriptstyle \left(A,H,\bsmat {\scriptscriptstyle C}\\[0.3mm]{\scriptscriptstyle E} \esmat,\bsmat {\scriptscriptstyle G_y} \\ {\scriptscriptstyle G_z}\esmat\right)}+\gS^\star_{\scriptscriptstyle (A,H,C,G_y)} \Big), \label{neweq}
\eea
where we used Theorem~\ref{the14} and, again, Theorem~\ref{the13d}.
From 
\[
\gV^\star_{\scriptscriptstyle \left(A,H,\bsmat {\scriptscriptstyle C}\\[0.3mm]{\scriptscriptstyle E} \esmat,\bsmat {\scriptscriptstyle G_y} \\ {\scriptscriptstyle G_z}\esmat\right)}\subseteq 
\gV^\star_{\scriptscriptstyle \left(A,H,E,G_z\right)}\subseteq \gV^\star_{\scriptscriptstyle \left(A,[\,B\;\;H\,],E, [\,D_z\;\;G_z\,]\right)}=\gV^\star_{\scriptscriptstyle (A,B,E,D_z)},
\]
where the equality comes from Theorem~\ref{the14}, and
since $\gS^\star_{\scriptscriptstyle (A,H,C,G_y)}\subseteq \gV^\star_{\scriptscriptstyle (A,B,E,D_z)}$, we find $\gV^\star_{\scriptscriptstyle \left(A,H,\bsmat {\scriptscriptstyle C}\\[0.3mm]{\scriptscriptstyle E} \esmat,\bsmat {\scriptscriptstyle G_y} \\ {\scriptscriptstyle G_z}\esmat\right)}+\gS^\star_{\scriptscriptstyle (A,H,C,G_y)}\subseteq \gV^\star_{\scriptscriptstyle (A,B,E,D_z)}$. We use this in 
(\ref{neweq}) and obtain
\bea
\gV_m\cap \gS_M=\gS^\star_{\scriptscriptstyle (A,[\,B\;\;H\,],E,[\,D_z\;\;G_z\,]}
\cap
(\gV^\star_{\scriptscriptstyle \left(A,H,\bsmat {\scriptscriptstyle C}\\[0.3mm]{\scriptscriptstyle E} \esmat,\bsmat {\scriptscriptstyle G_y} \\ {\scriptscriptstyle G_z}\esmat\right)}+\gS^\star_{\scriptscriptstyle (A,H,C,G_y)}). \label{pippo}
\eea
Consider the other inclusion (together with Theorem~\ref{the13d})
$\gS^\star_{\scriptscriptstyle (A,H,C,G_y)}=
\gS^\star_{\scriptscriptstyle \left(A,H,\bsmat {\scriptscriptstyle C}\\[0.3mm]{\scriptscriptstyle E} \esmat,\bsmat {\scriptscriptstyle G_y} \\ {\scriptscriptstyle G_z}\esmat\right)}\subseteq 
\gS^\star_{\scriptscriptstyle \left(A,H,E,G_z\right)}\subseteq \gS^\star_{\scriptscriptstyle \left(A,[\,B\;\;H\,],E, [\,D_z\;\;G_z\,]\right)}$.
We can use the modular rule on (\ref{pippo}) to obtain
\bea
\label{pluto}
\gV_m\cap \gS_M=\gS^\star_{\scriptscriptstyle (A,H,C,G_y)} +(\gS^\star_{\scriptscriptstyle (A,[\,B\;\;H\,],E,[\,D_z\;\;G_z\,]}\cap\gV^\star_{\scriptscriptstyle \left(A,H,\bsmat {\scriptscriptstyle C}\\[0.3mm]{\scriptscriptstyle E} \esmat,\bsmat {\scriptscriptstyle G_y} \\ {\scriptscriptstyle G_z}\esmat\right)}),
\eea
where $\gS^\star_{\scriptscriptstyle (A,[\,B\;\;H\,],E,[\,D_z\;\;G_z\,]}\cap\gV^\star_{\scriptscriptstyle \left(A,H,\bsmat {\scriptscriptstyle C}\\[0.3mm]{\scriptscriptstyle E} \esmat,\bsmat {\scriptscriptstyle G_y} \\ {\scriptscriptstyle G_z}\esmat\right)}\subseteq C^{-1}\ima G_y$.
Adding $C^{-1}\ima G_y$ to both sides of (\ref{pluto}) yields
$(\gV_m\cap \gS_M)+C^{-1}\ima G_y=\gS^\star_{\scriptscriptstyle (A,H,C,G_y)}+C^{-1}\ima G_y$.
Thus,
\beann
\ima [\begin{array}{cccc} T_1 & T_3 \end{array}]&=&\gS_M\subseteq C^{-1}\ima G_y +  \gS^\star_{\scriptscriptstyle \left(A,H,C,G_y\right)}\\
&=&C^{-1}\ima G_y+(\gV_m\cap \gS_M)=C^{-1}\ima G_y+\ima T_1,
\eeann
so that it is always possible to choose $T_3$ in such a way that $\ima T_3\subseteq C^{-1}\ima G_y$.

Recall that $\gV_m+\gS_M$ is $(A,B,E,D_z)$-output nulling, see Lemma~\ref{lem22}; if we denote by $\gR_{\gV_m+\gS_M}^\star$ the output nulling reachability subspace on $\gV_m+\gS_M$, there holds $(\gS_M+\gV_m)\cap B\,\ker D_z=
\gR_{\gV_m+\gS_M}^\star\cap B\,\ker D_z$. Again, since $\gV_m+\gS_M$ is $(A,B,E,D_z)$-output nulling, the subspace $\gR_{\gV_m+\gS_M}^\star$ is contained in $\gR^\star_{\scriptscriptstyle (A,B,E,D_z)}$, which 
in turn is contained in $\gV_m$. Thus,
$\ima [\begin{array}{ccc} T_1 & T_2 & T_3 \end{array}]\cap B\,\ker D_z = \ima [\begin{array}{ccc} T_1 & T_2\end{array}]\cap B\,\ker D_z$.
Then, we can also choose $T_4$ so that $\ima [\begin{array}{ccc} T_1 & T_2 & T_4 \end{array}]\supseteq B\,\ker D_z$.
%
%
Let $A_1=T^{-1}\,A\,T$, $B_1=T^{-1} \,B$, $H_1=T^{-1}\,H$, $C_1=C\,T$, $E_1=E\,T$.
Taking $F_1$ to be an $(A,B,E,D_z)$-output nulling friend 
of $\gS_M+\gV_m$, we obtain 
{
\[
A_1+B_1\,F_1=
\bmat{cccc}
\star& \star & \star & \star \\[-1.5mm]
\star& \star & \star & \star \\[-1.5mm]
%
0 & 0& A_{\scriptscriptstyle 33} & A_{\scriptscriptstyle 34}+B_{\scriptscriptstyle 32}\,F_{\scriptscriptstyle 24}  \\[-1.5mm]
0  & 0  & 0  &\star \emat.
\]
}
Since $\gV_m\subseteq \bar{\gV}\subseteq \gS_M+\gV_m$, we can write $\bar{\gV}=\ima [\begin{array}{cc} T_1 & T_2 \end{array}]+T_3\,X$ for a certain matrix $X$.
In the new basis, we can write
\[
\bar{\gV}=\ima \bmat{ccc} I & 0 & 0 \\
0 & I & 0 \\ 0 & 0 & X \\ 0 & 0 & 0 \emat.
\]
 Since $\bar{\gV}$ is $(A_1+B_1\,F_1)$-invariant, there exists a matrix $M$ partitioned comformably such that 
 {
\[
\bmat{cccc}
 \star& \star & \star & \star \\[-1.5mm]
 \star& \star & \star & \star \\[-1.5mm]
 0 & 0& A_{\scriptscriptstyle 33} & \star \\[-1.5mm]
 0  & 0  & 0  & A_{\scriptscriptstyle 44}+B_{\scriptscriptstyle 41} F_{\scriptscriptstyle 14}+B_{\scriptscriptstyle 42} F_{\scriptscriptstyle 24}  \emat \bmat{ccc}    I & 0 & 0     \\[-1.5mm]
    0 & I & 0     \\[-1.5mm]     0 & 0 & X     \\[-1.5mm]     0 & 0 & 0     \emat= \bmat{ccc}     I & 0 & 0     \\[-1.5mm]
    0 & I & 0     \\[-1.5mm]     0 & 0 & X     \\[-1.5mm]     0 & 0 & 0     \emat \underbrace{\bmat{ccc}  M_{\scriptscriptstyle 11} & M_{\scriptscriptstyle 12} & M_{\scriptscriptstyle 13}    \\[-1.5mm]
   M_{\scriptscriptstyle 21} & M_{\scriptscriptstyle 22} & M_{\scriptscriptstyle 23}   \\[-1.5mm]
  M_{\scriptscriptstyle 31} & M_{\scriptscriptstyle 32} & M_{\scriptscriptstyle 33}   \emat}_{M},
\]
}
from which we find $A_{33}\,X=X\,M_{33}$. Hence, $\ima X$ is $A_{3,3}$-invariant, and since $\bar{\gV}$ is internally stabilizable, then $\ima X$ is an internally stable $A_{3,3}$-invariant.

Similarly, choosing a friend $G$ of $\gV_m\cap \gS_M$ we obtain
{
 \[
A_1+G_1\,C_1=\bmat{cccc}
 \star &  \star  & A_{\scriptscriptstyle 13}+G_{\scriptscriptstyle 11}\,C_{\scriptscriptstyle 13} & \star\\[-1mm]
0 & \star  & 0  & \star \\[-1mm]
0 &  \star  & A_{\scriptscriptstyle 33} & \star  \\[-1mm]
0 &  \star  & 0& \star \emat.
\]
}
Since $\gS_M\cap \gV_m \subseteq \bar{\gS}\subseteq \gS_M$, we can write $\bar{\gS}=\ima T_1 +T_3\,Y$ for a certain matrix $Y$. Since $\bar{\gS}\subseteq \bar{\gV}$, then $T_3\,Y \subseteq T_3\,X$; thus $\ima Y \subseteq \ima X$. In the new basis, we can write
\[
\bar{\gS}=\ima \bmat{cc} I & 0  \\
0 & 0 \\ 0 & Y \\ 0 & 0 \emat.
\]
From the $(A_1+G_1\,C_1)$-invariance of $\bar{\gS}$, there exists a matrix $N$ partitioned comformably such that \\[-0.2cm]
\[
\bmat{cccc}
 \star &  \star  & A_{\scriptscriptstyle 13}+G_{\scriptscriptstyle 11}\,C_{\scriptscriptstyle 13} & \star\\[-1.5mm]
0 & \star  & 0  & \star \\[-1.5mm]
0 &  \star  & A_{\scriptscriptstyle 33} & \star  \\[-1.5mm]
0 &  \star  & 0& \star \emat
\bmat{ccc} I & 0  \\[-1.5mm]
0 & 0 \\[-1.5mm] 0 & Y \\[-1.5mm] 0 & 0 \emat=
\bmat{ccc} I & 0  \\[-1.5mm]
0 & 0 \\[-1.5mm] 0 & Y \\[-1.5mm] 0 & 0 \emat\underbrace{\bmat{cc} N_{\scriptscriptstyle 11} & N_{\scriptscriptstyle 12} \\[-1.5mm] N_{\scriptscriptstyle 21} & N_{\scriptscriptstyle 22} \emat}_{N},
\]
which yields\ $A_{33}\,Y=Y\,N_{22}$. Hence, $\ima Y$ is $A_{33}$-invariant, and it is externally stabilizable since $\bar{\gS}$ is externally detectable. Now consider a further change of basis for $A_{33}$ given by $\tilde{T}=[\begin{array}{ccc} \tilde{T}_1 & \tilde{T}_2 & \tilde{T}_3 \end{array}]$, where $\tilde{T}_1$ is a basis for $\ima Y$, and $[\begin{array}{cc} \tilde{T}_1 & \tilde{T}_2 \end{array}]$ is a basis for $\ima X$. Then
\[
\tilde{T}^{-1}\,A_{33}\,\tilde{T}=\bmat{ccc} A_{33}^1 & \star & \star \\ 0 & A_{33}^2 & \star \\ 0 & 0 & A_{33}^3\emat.
\]
Since $\ima X$ is internally stabilizable, $A_{33}^1$ and $ A_{33}^2 $ are stable; Since $\ima Y$ is externally stabilizable, $A_{33}^2$ and $ A_{33}^3 $ are stable. It follows that $A_{33}$ is stable, so that $\gS_M+\gV_m$ is internally stabilizable.
\endproof

\section*{Concluding remarks}
In this paper, we have developed a geometric solution to the disturbance decoupling by dynamic output feedback for systems which are not necessarily strictly proper, using the notions of self boundedness and self hiddenness. The building blocks of this solution do not require eigenspace computations that are at the basis of a solution involving stabilizability and detectability subspaces: the solution given here remains in the realm of finite arithmetics. 
The crucial issue in the extension of the classical theory to the nonstrictly proper case is the well-posedness of the closed-loop, which has to be handled separately from the other solvability conditions. Importantly, in this paper we have showed that checking this condition for the pair of subspaces $\gV^\star_{\scriptscriptstyle (A,B,E,D_z),g}$ and $\gS^\star_{\scriptscriptstyle (A,H,C,G_y),g}$, or for the pair of subspaces $\gV_m+\gS_M$ and $\gS_M$, is equivalent to checking the same condition for the pair of subspaces $\gV^\star_{\scriptscriptstyle (A,B,E,D_z)}$ and $\gS^\star_{\scriptscriptstyle (A,H,C,G_y)}$.

\section*{Appendix A: projection and intersection}

Consider two vector spaces $\gX$ and $\gP$. Let $\gS$ be a subspace of
$\gX\oplus \gP$. The linear operators $\mathfrak{p},\mathfrak{i}$ 
 are defined as
\beann
\mathfrak{p}\,(\gS) \ns&\ns \defi \ns&\ns  \left\{\x\in \gX\,\Big|\,\exists \,\p  \in \gP\,:\; \bsmat \x \\[1mm] \p \esmat \in \gS\right\} \\
\mathfrak{i}\,(\gS)  \ns&\ns \defi \ns&\ns  \left\{\x\in \gX\,\Big|\, \bsmat \x \\[1mm] 0 \esmat \in \gS\right\},
\eeann
where $\mathfrak{p}\,(\gS) $ is referred to as the projection of $\gS$ on $\gX$ and  
 $\mathfrak{i}\,(\gS) $ is the intersection of $\gS$ with $\gX$. It is easy to see that $
\mathfrak{p}\,(\gS) $ and $\mathfrak{i}\,(\gS) $ are subspaces of $\gX$. Both operators preserve addition and intersection, and 
 $\mathfrak{p}\,(\gW^\perp) = \bigl(\mathfrak{i}\,(\gW)\bigr)^\perp$, see \cite[Prop.~5.1.3]{Basile-M-92}.

\begin{lemma}
\label{PI1}
Let $\gW\supseteq \ima \bsmat H_1 \\[1mm] H_2 \esmat$. Then, $\mathfrak{p}\,(\gW) \supseteq \ima H_1$.
\end{lemma}
\ \\[-1.4cm]
%

\begin{lemma}
\label{PI2}
Let $\gW\subseteq \ker [\begin{array}{cc} C_1 & C_2 \end{array}]$. Then, $\mathfrak{i}\,(\gW)\subseteq \ker C_1$.
\end{lemma}
%

%

\section*{Appendix B}
In this Appendix, we recall some fundamental geometric results for a quadruple $(A,B,C,D)$. These are restatements or dualizations of the results in \cite[Appx.~A]{DMCuellar-M-01} and \cite[Lemma 3]{N-TAC-08}, see also \cite[Sec.~5]{N-EJC-07}.
We begin by studying the inclusion $\ima L \subseteq \gV^\star_{\scriptscriptstyle (A,B,C,D)}$. 

\begin{theorem}{\bf \cite[Lem.~3]{N-TAC-08}}
\label{the11}
Let $\ima L \subseteq \gV^\star_{\scriptscriptstyle (A,B,C,D)}$. The following results hold:\\[-0.7cm]
\begin{description}
\item{{\em i)}} $\gV^\star_{\scriptscriptstyle (A,B,C,D)}={\gV}^\star_{\scriptscriptstyle (A,[\,B\;\;L\,],C,[\,D\;\; 0 \,])}$;
\item{{\em ii)}} $\Phi_{\scriptscriptstyle (A,[\,B\;\;L\,],C,[\,D\;\; 0 \,])}\subseteq \Phi_{\scriptscriptstyle (A,B,C,D)}$;
\item{{\em iii)}} $\ima L   \subseteq \gV$ $\quad \forall\,{\gV}\in\Phi(A,[\,B\;\;L\,],C,[\,D\;\; 0 \,])$.
\end{description}
\end{theorem}

\begin{theorem}{\bf \cite[Prop.~A.1]{DMCuellar-M-01}}
\label{the12}
$\ima L \subseteq {\gV}^\star_{\scriptscriptstyle (A,B,C,D)}$ if and only if  $\ima L \subseteq \gR^\star_{\scriptscriptstyle (A,[\,B\;\;L\,],C,[\,D\;\; 0 \,])}$.
\end{theorem}

\begin{theorem}{\bf \cite[Prop.~A.3]{DMCuellar-M-01}}
\label{the13}
If $\ima L \subseteq \gV^\star_{\scriptscriptstyle (A,B,C,D)}$,
the subspace $\gR^\star_{\scriptscriptstyle (A,[\,B\;\;L\,],C,[\,D\;\; 0 \,])}$ is the smallest of all the $(A,B,C,D)$-self bounded subspaces $\gV$ satisfying $\ima L \subseteq \gV$.
\end{theorem}

The following three results are a generalization of the last three: they are concerned with a geometric condition in the form $\ima \bsmat  L_1 \\[1mm] L_2 \esmat \subseteq (\gV^\star_{\scriptscriptstyle (A,B,C,D)} \oplus 0_{\scriptscriptstyle \gY})+\ima \bsmat B \\[1mm] D \esmat$ which arises in the solution of the decoupling of a measurable disturbance. 

\begin{theorem}
\label{the14}
Let $\ima \bsmat  L_1 \\[1mm] L_2 \esmat \subseteq (\gV^\star_{\scriptscriptstyle (A,B,C,D)} \oplus 0_{\scriptscriptstyle \gY})+\ima \bsmat B \\[1mm] D \esmat$.
The following results hold:\\[-0.7cm]

\begin{description}
\item{{\em i)}} $\gV^\star_{\scriptscriptstyle (A,B,C,D)}={\gV}^\star_{\scriptscriptstyle (A,[\,B\;\;L_1\,],C,[\,D\;\; L_2 \,])}$;
\item{{\em ii)}} $\Phi_{\scriptscriptstyle (A,[\,B\;\;L_1\,],C,[\,D\;\; L_2 \,])}\subseteq \Phi_{\scriptscriptstyle (A,B,C,D)}$;
\item{{\em iii)}} $\ima \bsmat L_1 \\[1mm] L_2 \esmat \subseteq \gV \oplus 0_{\scriptscriptstyle \gY}+\ima \bsmat B \\[1mm] D \esmat$  $\quad \forall\,{\gV}\in\Phi_{\scriptscriptstyle (A,[\,B\;\;L_1\,],C,[\,D\;\; L_2 \,])}$.
\end{description}
\end{theorem}

\begin{theorem}
$\ima \bsmat  L_1 \\[1mm] L_2 \esmat \subseteq (\gV^\star_{\scriptscriptstyle (A,B,C,D)} \oplus 0_{\scriptscriptstyle \gY})+\ima \bsmat B \\[1mm] D \esmat$
if and only if  
$\ima \bsmat L_1 \\[1mm] L_2 \esmat \subseteq (\gR^\star_{\scriptscriptstyle (A,[\,B\;\;L_1\,],C,[\,D\;\;L_2\,])} \oplus 0_{\scriptscriptstyle \gY})+\ima \bsmat B \\[1mm] D \esmat$.
\end{theorem}

\begin{theorem}
If
$\ima \bsmat  L_1 \\[1mm] L_2 \esmat \subseteq (\gV^\star_{\scriptscriptstyle (A,B,C,D)} \oplus 0_{\scriptscriptstyle \gY})+\ima \bsmat B \\[1mm] D \esmat$,
the subspace $\gR^\star_{\scriptscriptstyle (A,[\,B\;\;L_1\,],C,[\,D\;\; L_2 \,])}$ is the smallest of all the $(A,B,C,D)$-self bounded subspaces $\gV$ satisfying
$\ima \bsmat L_1 \\[1mm] L_2 \esmat \subseteq (\gV \oplus 0_{\scriptscriptstyle \gY})+\ima \bsmat B\\[1mm] D \esmat$.
\end{theorem}

We now dualize all the previous results. The first three involve an inclusion $\gS_{\scriptscriptstyle (A,B,C,D)}^\star\subseteq \ker M$, for some matrix $M$.

\begin{theorem}
\label{the11d}
Let $\gS_{\scriptscriptstyle (A,B,C,D)}^\star\subseteq \ker M$. The following hold:
\begin{description}
\item{{\em i)}} $\gS^\star_{\scriptscriptstyle (A,B,C,D)}={\gS}^\star_{\scriptscriptstyle \left(A,B,\bsmat {\scriptscriptstyle C} \\ {\scriptscriptstyle M} \esmat,\bsmat {\scriptscriptstyle D} \\ {\scriptscriptstyle 0} \esmat\right)}$;
\item{{\em ii)}} $\Psi_{\scriptscriptstyle \left(A,B,\bsmat {\scriptscriptstyle C}\\ {\scriptscriptstyle M} \esmat,\bsmat {\scriptscriptstyle D} \\ {\scriptscriptstyle 0} \esmat\right)} \subseteq \Psi_{\scriptscriptstyle (A,B,C,D)}$;
\item{{\em iii)}} $\gS   \subseteq \ker M$ $\quad \forall\,{\gS}\in\Psi_{\scriptscriptstyle \left(A,B,\bsmat {\scriptscriptstyle C}  \\{\scriptscriptstyle M} \esmat,\bsmat {\scriptscriptstyle D} \\ {\scriptscriptstyle 0} \esmat\right)}$.
\end{description}
\end{theorem}

\begin{theorem}
\label{the12d}
${\gS}^\star_{\scriptscriptstyle (A,B,C,D)}\subseteq \ker M$ iff  $\gQ^\star_{\scriptscriptstyle \left(A,B,\bsmat {\scriptscriptstyle C}\\ {\scriptscriptstyle M} \esmat,\bsmat {\scriptscriptstyle D} \\ {\scriptscriptstyle 0} \esmat\right)}\subseteq \ker M$.
\end{theorem}

\begin{theorem}
\label{the13d}
If $\gS^\star_{\scriptscriptstyle (A,B,C,D)} \subseteq \ker M$,
the subspace $\gQ^\star_{\scriptscriptstyle \left(A,B,\bsmat {\scriptscriptstyle C} \\ {\scriptscriptstyle M} \esmat,\bsmat {\scriptscriptstyle D} \\ {\scriptscriptstyle 0} \esmat\right)}$ is the largest of all the $(A,B,C,D)$-self hidden subspaces $\gS$ satisfying
$\gS\subseteq \ker M$.
\end{theorem}

Finally, we consider the generalization $(\gS_{\scriptscriptstyle (A,B,C,D)}^\star \oplus \gU) \cap \ker [\begin{array}{cc} C & D \end{array}]\subseteq \ker [\begin{array}{cc} M_1 & M_2 \end{array}]$ of the condition $\gS^\star_{\scriptscriptstyle (A,B,C,D)} \subseteq \ker M$. 

\begin{theorem}
\label{the11d1}
Let $(\gS_{\scriptscriptstyle (A,B,C,D)}^\star \oplus \gU) \cap \ker [\begin{array}{cc} C & D \end{array}]\subseteq \ker [\begin{array}{cc} M_1 & M_2 \end{array}]$. 
 The following results hold:
\begin{description}
\item{{\em i)}} $\gS^\star_{\scriptscriptstyle (A,B,C,D)}={\gS}^\star_{\scriptscriptstyle \left(A,B,\bsmat {\scriptscriptstyle C} \\ {\scriptscriptstyle M_2} \esmat,\bsmat {\scriptscriptstyle D} \\ {\scriptscriptstyle M_2} \esmat\right)}$;
\item{{\em ii)}} $\Psi_{\scriptscriptstyle \left(A,B,\bsmat {\scriptscriptstyle C}\\{\scriptscriptstyle M_1} \esmat,\bsmat {\scriptscriptstyle D} \\ {\scriptscriptstyle M_2} \esmat\right)} \subseteq \Psi_{\scriptscriptstyle (A,B,C,D)}$;
\item{{\em iii)}} $(\gS   \oplus \gU) \cap \ker [\begin{array}{cc} C & D \end{array}]\subseteq \ker [\begin{array}{cc} M_1 & M_2 \end{array}]$ $\forall\,{\gS}\in\Psi_{\scriptscriptstyle \left(A,B,\bsmat {\scriptscriptstyle C} \\{\scriptscriptstyle M_1} \esmat,\bsmat {\scriptscriptstyle D} \\ {\scriptscriptstyle M_2} \esmat\right)}$.
\end{description}
\end{theorem}

\begin{theorem}
$(\gS_{\scriptscriptstyle (A,B,C,D)}^\star \oplus \gU) \cap \ker [\begin{array}{cc} C & D \end{array}]\subseteq \ker [\begin{array}{cc} M_1 & M_2 \end{array}]$
 iff
$(\gQ^\star_{\scriptscriptstyle \left(A,B,\bsmat {\scriptscriptstyle C} \\ {\scriptscriptstyle M_1} \esmat,\bsmat {\scriptscriptstyle D} \\ {\scriptscriptstyle M_2} \esmat\right)} \oplus \gU) \cap \ker [\begin{array}{cc} C & D \end{array}]\subseteq \ker [\begin{array}{cc} M_1 & M_2 \end{array}]$.
\end{theorem}

\begin{theorem}
\label{the12d3}
If $(\gS_{\scriptscriptstyle (A,B,C,D)}^\star \oplus \gU) \cap \ker [\begin{array}{cc} C & D \end{array}]\subseteq \ker [\begin{array}{cc} M_1 & M_2 \end{array}]$, the subspace $\gQ^\star_{\scriptscriptstyle \left(A,B,\bsmat {\scriptscriptstyle C}\\ {\scriptscriptstyle M_1} \esmat,\bsmat {\scriptscriptstyle D} \\ {\scriptscriptstyle M_2} \esmat\right)}$ is the largest of all the $(A,B,C,D)$-self hidden subspaces $\gS$ satisfying $(\gS \oplus \gU) \cap \ker [\begin{array}{cc} C & D \end{array}]\subseteq \ker [\begin{array}{cc} M_1 & M_2 \end{array}]$.
\end{theorem}

\end{document}